\documentclass[graybox]{svmult}

\usepackage{mathptmx}       
\usepackage{helvet}         
\usepackage{courier}        
\usepackage{type1cm}        
\usepackage{makeidx}         
\usepackage{graphicx}        
\usepackage{multicol}        
\usepackage[bottom]{footmisc}

\usepackage{amsmath}
\usepackage{amssymb}
\usepackage{amscd}
\usepackage{indentfirst}

\def\1{{\bf 1}}

\makeindex             

\begin{document}

\title*{Recent developments of Lyapunov-type inequalities for fractional differential equations}
\titlerunning{Fractional Lyapunov-type inequalities}
\author{Sotiris K. Ntouyas, Bashir Ahmad and Theodoros P. Horikis}
\authorrunning{S.K. ~Ntouyas, B.~Ahmad and T.P.~Horikis }
\institute{S.K. ~Ntouyas  \at Department of Mathematics, University of  Ioannina, 451 10 Ioannina Greece\\ \email{sntouyas@uoi.gr}\\
B.~Ahmad \at Nonlinear Analysis and Applied Mathematics (NAAM)-Research Group,
Department of Mathematics, Faculty of Science, King Abdulaziz
University,\\  P.O. Box 80203, Jeddah 21589, Saudi Arabia\\ \email{bashirahmad$_{-}$qau@yahoo.com}\\
T.P.~Horikis \at Department of Mathematics, University of  Ioannina, 451 10 Ioannina Greece\\ \email{horikis@uoi.gr}}

\maketitle

\abstract{A survey of results on Lyapunov-type  inequalities for   fractional differential equations associated with a variety of boundary conditions is presented. This includes Dirichlet, mixed, Robin, fractional, Sturm-Liouville, integral, nonlocal, multi-point, anti-periodic, conjugate, right-focal and impulsive conditions.  Furthermore, our study includes   Riemann-Liouville, Caputo, Hadamard, Prabhakar, Hilfer and conformable type fractional derivatives. Results for boundary value problems involving fractional $p$-Laplacian, fractional operators with nonsingular Mittag-Leffler kernels, $q$-difference, discrete, and impulsive equations, are also taken into account.}


\section{Introduction and Preliminaries}
Integral inequalities are fundamental  in the
study of  quantitative properties of solutions of differential and integral equations. The
Lyapunov-type inequality is one of such inequalities when investigating the zeros of solutions of differential equations.
A method for deriving a Lyapunov-type inequality for  boundary value problems dates back to Nehari \cite{Ne} and is based on the idea of converting the given problem into an integral equation. To illustrate this method, let us consider the following    boundary value problem:
 \begin{equation}\label{E-1}
 \left\{\begin{array}{ll}
 y''(t)+q(t)y(t)=0, \quad a<t<b,\\
 y(a)=y(b)=0,
 \end{array}
 \right.
 \end{equation}
 where $a,b\in {\mathbb R}, a<b$ are consecutive zeros of $y(t)$ and $y(t)\not\equiv 0,$ for all $t\in (a,b).$  It can easily be shown  that    problem (\ref{E-1}) is equivalent to the integral equation
\begin{equation}\label{E-2}
y(t)=\int_a^b G(t,s) q(s) y(s) ds,
\end{equation}
where $G(t,s)$ is the Green's function given by
\begin{equation}\label{green-1}
G(t,s)=-
\begin{cases}
\displaystyle\frac{(t-a)(b-s)}{b-a}, &a\leq s\leq t\leq b,\\
\displaystyle\frac{(s-a)(b-t)}{b-a}, &a\leq t \leq s \leq b.
\end{cases}
\end{equation}
 Taking the absolute value of both sides of equation (\ref{E-2}), and taking into account that $y$ does not have any zeros in $(a,b),$we get
\begin{equation}\label{E-3}
1\le \max_{a\le t\le b}\int_a^b |G(t,s)| |q(s)|  ds,
\end{equation}
which yields the desired Lyapunov inequality
\begin{equation}\label{E-4}
 \int_a^b  |q(s)|  ds\ge \frac{1}{\max_{(t,s)\in [a,b]\times [a,b]}|G(t,s)|}.
\end{equation}
In the special case  where  can find $H(s)$ explicily such that $\displaystyle \max_{t\in [a,b]}|G(t,s)|\le H(s),$ then we obtain the following inequality
$$ 1\le \int_a^b H(s)|q(s)|ds.$$
Clearly the function $H(t)$ for problem (\ref{E-1}) is $\displaystyle\frac{(t-a)(b-t)}{b-a}.$ Moreover, if we take the absolute maximum of the function $H(t)$ for all $t\in [a,b],$ then it is obtained the following well-known Lyapunov    inequality   \cite{Lya}.

\begin{theorem}\label{T-L} If  the boundary-value problem (\ref{E-1})
has a nontrivial solution, where $q$ is a real and continuous function, then
 \begin{equation}\label{E-5}
 \int_a^b|q(s)|ds>\frac{4}{b-a}.
  \end{equation}
\end{theorem}
The factor $4$ in the above inequality is sharp
and cannot be replaced by a larger number.

 Later Wintner in \cite{Win} and more authors thereafter generalized this result by replacing
  the function $|q(t)|$ in (\ref{E-5}) by the function $q^+(t),$  $q^+(t)=\max\{q(t),0\},$ where now the resulting inequality reads:
\begin{equation}\label{E-5a}
\int_a^bq^{+}(s)ds>\frac{4}{b-a},
\end{equation}
with $q^{+}(t)=\max\{q(t),0\}.$

In \cite{Har}, Hartman expanded further this result with   the following inequality
\begin{equation}\label{E-5b}
\int_a^b  (b-t)(t-a)q^{+}(t)dt>(b-a),
\end{equation}
which is sharper than both  (\ref{E-5}) and  (\ref{E-5a}).

Clearly, (\ref{E-5b})  implies (\ref{E-5a}) as $\displaystyle  (b-t)(t-a)\le \frac{(b-a)^2}{4}$ for all $t\in [a,b]$  and the equality holds when $\displaystyle t=\frac{a+b}{2}.$

 It is worth mentioning   that  inequality (\ref{E-5}) has found many practical applications in differential equations (oscillation theory, disconjugacy, eigenvalue problems, etc.), for instance, see \cite{Ca}-\cite{ZhTa} and references therein.  A thorough literature review  dealing with  continuous
and discrete Lyapunov inequalities and their applications can be found in  \cite{Ch} and \cite{BH} (which also includes an excellent account on the history of such inequalities).

In many engineering and scientific disciplines such as physics, chemistry, aerodynamics,
electrodynamics of complex media, polymer rheology, economics, control theory,
signal and image processing, biophysics, blood flow and related phenomena,  fractional differential
and integral equations represent  processes in a more effective manner than
 their integer-order counterparts. This aspect has led to the  increasing popularity in the study  of fractional order differential and integral equations among mathematicians and researchers.
In view of their extensive applications in various fields,  the topic of inequalities for fractional differential equations  has also attracted a significant attention in recent years.

This survey article  is organized as follows. In Section 2 we introduce the reader to  some basic concepts of fractional calculus. In Section 3 we summarize Lyapunov-type inequalities for fractional boundary value problems with different  kinds of boundary  conditions. In Section 4 we consider the inequalities for nonlocal and multi-point boundary value problems. Results on $p$-Laplacians are discussed in Section 5, while results on mixed fractional derivatives are given in Section 6. Section 7 deals with Lyapunov-type inequalities for Hadamard fractional differential equations. In Section 8,  inequalities involving Prabhakar fractional differential equations  are discussed. Section 9 contains the results on  fractional $q$-difference equations, while Section 10 consists of the results involving fractional derivatives with respect to a certain function. Inequalities involving left and right derivatives, operators with nonsingular  Mittag-Leffler kernels,  discrete fractional  differential equations,  and impulsive fractional boundary value problems are respectively given in Sections 11, 12, 13 and 14 respectively. We include the results for Hilfer rand Katugampola fractional differential equations in Sections 15 and 16 respectively, and  conclude with Section 17 with results on  conformable fractional differential equations. Note that our goal here is a more complete and comprehensive review and as such the choice is made to include as many results as possible to illustrate the progress on the matter. Any proofs (which are rather long) are omitted, for this matter, and the reader is referred to the relative article accordingly.

\section{Fractional Calculus}
Here we introduce some basic definitions of fractional calculus \cite{Kil, pod} and recall some results that we need in the sequel.

\begin{definition} (Riemann-Liouville fractional integral)
Let $\alpha\ge 0$ and $f$ be a real function defined on $[a, b].$ The Riemann-Liouville fractional integral of order $\alpha$ is defined by $(I^0 f )(x) = f (x)$ and
$$(I^{\alpha} f)(t)=\frac{1}{\Gamma(\alpha)}\int_a^t(t-s)^{\alpha-1}f(s)ds, ~~ \alpha>0, ~~ t\in [a,b]$$
provided the right hand side is point-wise defined on $[0,\infty),$ where $\Gamma(\alpha)$ is the Euler Gamma function: $\displaystyle \Gamma(\alpha)=\int_0^{\infty}t^{\alpha-1}e^{-t}dt.$
\end{definition}
\begin{definition} (Riemann-Liouville fractional derivative)
The Riemann-Liouville fractional derivative of order $\alpha\ge 0$ is defined
by $(D^0 f)(t) = f (t)$ and
$$(D^{\alpha} f)(t) = (D^m\, I^{m-\alpha} f)(t)~\mbox{for}~ \alpha > 0,$$
where $m$ is the smallest integer greater than or equal  to $\alpha.$
\end{definition}
\begin{definition} (Caputo fractional derivative)
The Caputo fractional derivative of order $\alpha\ge 0$ is defined
by $(^CD^0 f)(t) = f (t)$ and
$$(^CD^{\alpha} f)(t) = (I^{m-\alpha} D^m\, f)(t)~\mbox{for}~ \alpha > 0,$$
where $m$ is the smallest integer greater than or equal  to $\alpha.$
\end{definition}

Notice that the differential operators of arbitrary order are nonlocal in nature  and appear in the mathematical modeling of several real world phenomena due to this characteristic (see e.g. \cite{Kil}).

 If $f \in C([a, b], {\mathbb R}),$  then the Riemann-Liouville fractional integral of order $\gamma>0$  exists on $[a, b].$ On the other hand, following \cite[Lemma 2.2, p. 73]{Kil}, we know that the Riemann-Liouville fractional derivative of
order $\gamma \in [n-1, n)$ exists almost everywhere (a.e.) on $[a, b]$ if $f \in AC^n([a, b], {\mathbb R}),$ where $C^k([a, b], {\mathbb R})$ $(k = 0, 1, \ldots)$ denotes the set of
 $k$ times continuously differentiable mappings on $[a, b],$ $AC([a, b], {\mathbb R})$ is the space of functions which are absolutely
continuous on $[a, b]$ and $AC^{(k)}([a, b], {\mathbb R})$ $ (k = 1, \ldots)$ is the space of functions $f$ such that $f \in C^{k-1}([a, b], {\mathbb R})$ and
$f^{(k-1)} \in AC([a, b], {\mathbb R}).$  In particular, $AC([a, b], {\mathbb R}) = AC^1([a, b], {\mathbb R}).$ (We recall here that  $AC([a,b], {\mathbb R})$ is the space of of functions $f$ which are absolutely continuous on $[a,b],$ and $ AC^n([a, b], {\mathbb R})$ the space of functions $f$ which have continuous derivatives up to order $n-1$ on $[a,b]$ such that $f^{(n-1)}(t)\in  AC([a, b], {\mathbb R})$).

\vskip 0.3cm

Now we enlist some important results involving fractional order operators \cite{Kil}.
\begin{proposition}
Let $f$ be a continuous function on some interval $J$ and $p, q>0.$ Then
$$(I^p~ I^q f)(t)=(I^{p+q}f)(t)=(I^q~ I^pf)(t)~ \mbox{on}~ J.$$
\end{proposition}

\begin{proposition}
Let $f$ be a continuous function on some interval $I$ and $\alpha\ge 0.$ Then
$$({\mathcal D}^{\alpha} I^{\alpha}f)(t)=f(t)~ \mbox{on}~ I,$$
with $\mathcal{D}$ being the Riemann-Liouville or Caputo fractional derivative operator.
\end{proposition}

\begin{proposition}
The general solution  of the following fractional differential equation
$$(D^qy)(t)=f(t), ~ t>a, ~ 0<q\le 1,$$
is $y(t)=c(t-a)^{q-1}+(I^qf)(t),~ c\in {\mathbb R}.$
\end{proposition}

\begin{proposition}
The general solution  of the following fractional differential equation
$$(^CD^qy)(t)=f(t), ~ t>a, ~ 0<q\le 1,$$
is $y(t)=c+(I^qf)(t),~ c\in {\mathbb R}.$
\end{proposition}

\section{Lyapunov-type inequalities for fractional differential equations with different boundary conditions}

Lyapunov-type inequalities involving
fractional differential operators have been investigated by many researchers in the recent years.   In 2013, Ferreira \cite{Fe} derived a Lyapunov-type inequality for   Riemann-Liouville fractional   differential equation with {\em Dirichlet boundary conditions}:
  \begin{equation}\label{L-0}
 \left\{\begin{array}{ll}
D^{\alpha}y(t)+q(t)y(t)=0, \quad a<t<b,\\[0.2cm]
 y(a)=y(b)=0,
 \end{array}
 \right.
 \end{equation}
 where $D^{\alpha}$ is the Riemann-Liouville fractional derivative of order  $1<\alpha\le 2$ and $q: [a,b]\to {\mathbb R}$ is a continuous function.

An appropriate approach for obtaining the Lyapunov inequality within  the framework of fractional differential equations relies on the idea of converting the boundary value problem  into an equivalent integral equation and then finding  the maximum value of its kernel function (Green's function).

It is straightforword  to show that the boundary value problem (\ref{L-0})
is equivalent to the integral equation
\begin{equation}\label{G-L-0}
y(t)=\int_a^bG(t,s)q(s)y(s)ds,
\end{equation}
where $G(t,s)$ is the Green's function defined by
\begin{equation}\label{Gr-L-0}
G(t,s)=\frac{1}{\Gamma(\alpha)}
\begin{cases}
\displaystyle\frac{(b-s)^{\alpha-1}(t-a)^{\alpha-1}}{(b-a)^{\alpha-1}}, &a\leq t\leq s\leq b,\\[0.3cm]
\displaystyle\frac{(b-s)^{\alpha-1}(t-a)^{\alpha-1}}{(b-a)^{\alpha-1}}-(t-s)^{\alpha-1}, &a\leq s \leq t \leq b.
\end{cases}
\end{equation}

Observe that the  Green's function  \eqref{Gr-L-0}   satisfies the following properties:
\begin{itemize}
\item[(i)] $G(t,s)\ge 0$, $\forall t,s\in [a,b];$
\item[(ii)]  $\displaystyle  \max_{s\in[a,b]} G(t,s)=G(s,s), s\in [a,b];$
\item[(iii)]~ $G(s,s)$ has a unique maximum, given by
$$\displaystyle  \max_{s\in[a,b]} G(s,s) =G\Big(\frac{a+b}{2},\frac{a+b}{2}\Big)=\frac{1}{\Gamma(\alpha)}\Big(\frac{b-a}{4}\Big)^{\alpha-1}.$$
\end{itemize}

 The Lyapunov inequality for  problem (\ref{L-0}) can be expressed as follows.

 \begin{theorem}\label{L=0}
If $y$ is a nontrivial solution of  the
 boundary value problem (\ref{L-0}), then
 \begin{equation}\label{L0}
 \int_a^b|q(s)|ds>\Gamma(\alpha)\Big(\frac{4}{b-a}\Big)^{\alpha-1}. \end{equation}
 \end{theorem}
Here we remark that Lyapunov's standard inequality (\ref{E-5}) follows   by taking $\alpha=2$ in the above inequality. Also, inequality (\ref{L0}) can be used to determine
intervals for the real zeros of the Mittag-Leffler function:
$$E_{\alpha}(z)=\sum_{k=1}^{\infty}\frac{z^k}{\Gamma(k\alpha+\alpha)}, ~ z\in {\mathbb C}, ~ \mathfrak{R}(\alpha)>0.$$

For the sake of convenience, let us consider the following
fractional Sturm-Liouville eigenvalue problem (with  $a = 0$ and $b = 1$):
 \begin{equation}\label{ML-0}
 \left\{\begin{array}{ll}
^CD^{\alpha}y(t)+\lambda y(t)=0, \quad 0<t<1,\\
 y(0)=y(1)=0.
 \end{array}
 \right.
 \end{equation}
 By Theorem \ref{L=0}, if $\lambda\in {\mathbb R}$ is an eigenvalue of  (\ref{ML-0}), that is, if $\lambda$ is a zero
of equation $E_{\alpha}(-\lambda)=0,$ then $|\lambda|>\Gamma(\alpha)4^{\alpha-1}.$ Therefore the Mittag-Leffler function
 $E_{\alpha}(z)$ has no real zeros for $|z|\le \Gamma(\alpha)4^{\alpha-1}.$

\vskip 0.3cm
 In 2014, Ferreira   \cite{Fe-1} replaced  the Riemann-Liouville fractional derivative in  problem (\ref{L-0})
 with Caputo  fractional derivative $^CD^{\alpha}$ and derived the following Lyapunov-type inequality for the resulting problem:
   \begin{equation}\label{L0a}
  \int_a^b|q(s)|ds>\frac{\Gamma(\alpha)\alpha^{\alpha}}{[(\alpha-1)(b-a)]^{\alpha-1}}. \end{equation}

\vskip 0.3cm
In 2015, Jleli and Samet \cite{JlSa} considered the fractional differential equation
\begin{equation}\label{eq7}
^CD^{\alpha}y(t)+q(t)y(t)=0, ~~~ 1<\alpha\le 2, ~~~ a<t<b,
 \end{equation}
 associated with the {\em mixed boundary conditions}
 \begin{equation}\label{eq8}
y(a)=0=y'(b)
 \end{equation}
 or
 \begin{equation}\label{eq9}
y'(a)=0=y(b).
 \end{equation}
An integral equation equivalent to   problem (\ref{eq7})-(\ref{eq8}) is
\begin{equation}\label{G-eq7}
y(t)=\int_a^bG(t,s)q(s)y(s)ds,
\end{equation}
where $G(t,s)$ is again the Green's function defined by
$$G(t,s)=\frac{H(t,s)}{\Gamma(\alpha)(b-s)^{2-\alpha}},$$
and
\begin{equation}\label{Gr-eq7}
H(t,s)=
\begin{cases}
\displaystyle (\alpha-1)(t-\alpha), &a\leq t\leq s\leq b,\\
\displaystyle  (\alpha-1)(t-\alpha)-(t-s)^{\alpha-1}(b-s)^{2-\alpha}, &a\leq s \leq t \leq b.
\end{cases}
\end{equation}
 The function $H$ satisfies the following inequality:
 $$|H(t,s)|\le \max\{(2-\alpha)(b-s), (\alpha-1)(s-a)\}~~\mbox{for all}~~ (t,s)\in [a,b]\times [a,b].$$

In relation to  problem (\ref{eq7})-(\ref{eq8}),  we have the following Lyapunov-type inequality.
\begin{theorem}\label{L-eq7}
If $y$ is a nontrivial solution of  the
 boundary value problem (\ref{eq7})-(\ref{eq8}), then
\begin{equation}\label{eq10}
 \int_a^b(b-s)^{\alpha-2}|q(s)|ds>\frac{\Gamma(\alpha)}{\max\{\alpha-1, 2-\alpha\}(b-a)}.
  \end{equation}
 \end{theorem}
In a similar manner, the Lyapunov-type inequality obtained for the boundary value problem (\ref{eq7})-(\ref{eq9}) is
   \begin{equation}\label{eq11}
 \int_a^b(b-s)^{\alpha-1}|q(s)|ds> \Gamma(\alpha).
  \end{equation}

As an application of  Lyapunov-type inequalities (\ref{eq10}) and (\ref{eq11}), we can obtain the intervals, where
certain Mittag-Leffler functions have no real zeros.
\begin{corollary}
Let $1<\nu\le 2.$ Then the Mittag-Leffler function $E_{\alpha}(z)$  has no real zeros for
$$z\in \Bigg(-\Gamma(\alpha)\frac{(\alpha-1)}{\max\{\alpha-1, 2-\alpha\}},0\Bigg].$$
\end{corollary}
The proof of the above corollary follows by applying  inequality (\ref{eq10}) to the following eigenvalue problem
 \begin{equation}\label{exam-10}
 \left\{\begin{array}{ll}
 (^CD^{\alpha} y)(t)+\lambda y(t)=0,  \quad 0<t<1,   \\[0.2cm]
\displaystyle y(0)=y'(1)=0.
 \end{array}
 \right.
 \end{equation}
 Moreover, by applying  inequality (\ref{eq11}) to the following eigenvalue problem
 \begin{equation}\label{exam-11}
 \left\{\begin{array}{ll}
 (^CD^{\alpha} y)(t)+\lambda y(t)=0,  \quad 0<t<1,   \\[0.2cm]
\displaystyle y'(0)=y(1)=0,
 \end{array}
 \right.
 \end{equation}
 we can obtain the following result:
 \begin{corollary}
Let $1<\nu\le 2.$ Then the Mittag-Leffler function $E_{\alpha}(z)$  has no real zeros for
$$z\in(\alpha\Gamma(\alpha),0].$$
\end{corollary}

In 2015, Jleli {\em et al.} \cite{L1}  obtained a Lyapunov-type inequality:
 {\small \begin{equation}\label{L1}
 \int_a^b(b-s)^{\alpha-2}(b-s+\alpha-1)|q(s)|ds\ge \frac{(b-a+2)\Gamma(\alpha)}{\rm{\max}\{b-a+1,((2-\alpha)/(\alpha-1))(b-a)-1\}},
 \end{equation}}
for the following problem with {\em Robin boundary conditions}
 \begin{equation}\label{L-1}
 \left\{\begin{array}{ll}
^CD^{\alpha}y(t)+q(t)y(t)=0, \,  1<\alpha\le 2, \quad a<t<b,\\[0.2cm]
 y(a)-y'(a)=y(b)+y'(b)=0,
 \end{array}
 \right.
 \end{equation}
 where $q: [a,b]\to {\mathbb R}$ is a continuous function.

 Using the  Lyapunov-type inequality (\ref{L1}), we can find an interval, where a linear combination of Mittag-Leffler functions  $\displaystyle E_{\alpha,\beta}=\sum_{k=0}^{\infty}\frac{z^k}{\Gamma(k\alpha+\beta)}, ~ \alpha>0, \beta>0, ~ z\in {\mathbb C}$ has no real zeros. In precise terms, we have the following result.
 \begin{theorem}
 Let $a<\alpha\le 2.$ Then $E_{\alpha,2}(z)+E_{\alpha,1}(z)+zE_{\alpha,\alpha}(z)$ has no real zeros for
 $$z\in \Bigg(\frac{-3\alpha\Gamma(\alpha)}{(1+\alpha)\max\{2, ((2-\alpha)/(\alpha-1)-1\}},0\Bigg].$$
 \end{theorem}

In 2015, Rong and Bai \cite{L3} considered a boundary value problems with {\em fractional boundary conditions}:
 \begin{equation}\label{L-3}
 \left\{\begin{array}{ll}
^CD^{\alpha}y(t)+q(t)y(t)=0, \, 1<\alpha\le 2, \quad a<t<b,\\[0.2cm]
\displaystyle y(a)=0, ~~ ^CD^{\beta}y(b)=0, \, 0<\beta\le 1,
 \end{array}
 \right.
 \end{equation}
 where  $q: [a,b]\to {\mathbb R}$ is a continuous function. The Lyapunov-type inequality derived for  problem (\ref{L-3}) is
 {\small \begin{equation}\label{L3}
 \int_a^b(b-s)^{\alpha-\beta-1}|q(s)|ds> \frac{(b-a)^{-\beta}}{\displaystyle\max\Big\{\frac{1}{\Gamma(\alpha)}-\frac{\Gamma(2-\beta)}{\Gamma(\alpha-\beta)}, \frac{\Gamma(2-\beta)}{\Gamma(\alpha-\beta)}, \frac{2-\alpha}{\alpha-1}\cdot \frac{\Gamma(2-\beta)}{\Gamma(\alpha-\beta)}\Big\}}.
 \end{equation}}
Later, Jleli and Samet \cite{L4}  obtained a Lyapunov-type inequality for a boundary value problems with {\em Sturm-Liouville boundary conditions}
 \begin{equation}\label{L-4}
 \left\{\begin{array}{ll}
^CD^{\alpha}y(t)+q(t)y(t)=0, \,  1< \alpha \le 2, \quad a<t<b,\\[0.2cm]
\displaystyle p y(a)-r y'(a)=y(b)=0, \, p>0, r\ge 0, \end{array}
 \right.
 \end{equation}
 where $q: [a,b]\to {\mathbb R}$ is a continuous function. The integral equation equivalent to the problem (\ref{L-4}) is
\begin{equation}\label{G-L-4}
y(t)=\int_a^bG(t,s)q(s)y(s)ds,
\end{equation}
where $G(t,s)$ is the Green's function defined by
\begin{equation}\label{Gr-L-4}
G(t,s)=\frac{1}{\Gamma(\alpha)}
\begin{cases}
\displaystyle\frac{(r/p+t-a)}{r/p+b-a}(b-s)^{\alpha-1}, &a\leq t\leq s\leq b,\\[0.2cm]
\displaystyle\frac{(r/p+t-a)}{r/p+b-a}(b-s)^{\alpha-1}-(t-s)^{\alpha-1}, &a\leq s \leq t \leq b.
\end{cases}
\end{equation}
In order to estimate this Green's function, we distinguish two cases:
  \begin{itemize}
 \item[(i)] If $\displaystyle \frac{r}{p}>\frac{b-a}{\alpha-1},$ then
 $0\le G(t,s)\le G(s,s), ~~(t,s)\in [a,b]\times [a,b]$ with
 $$\max_{s\in [a,b]}G(t,s)=\frac{1}{\Gamma(\alpha)}\frac{(r/p)(b-a)^{\alpha-1}}{\Big(r/p+b-a\Big)}.$$
  \item[(ii)]~ If $\displaystyle 0\le \frac{r}{p}\le\frac{b-a}{\alpha-1},$ then
  $\Gamma(\alpha)G(t,s)\le {\rm \max}\{{\mathcal A}(\alpha, r/p), {\mathcal B}(\alpha, r/p)\},$
   where
    $${\mathcal A}(\alpha, r/p)=\frac{(b-a)^{\alpha-1}}{(r/p+b-a)}\Bigg(\Big(\frac{(b-a)^{\alpha-1}}{(r/p+b-a)(\alpha-1)^{\alpha-1}}\Big)^{\frac{1}{\alpha-2}}(2-\alpha)-\frac{r}{p}\Bigg),$$
$$ {\mathcal B}(\alpha, r/p)=\Big(\frac{r}{p}+b-a\Big)^{\alpha-1}\frac{(\alpha-1){\alpha-1}}{\alpha^{\alpha}}.$$
 \end{itemize}
The Lyapunov inequalities corresponding to the above cases are given in the following result.
 \begin{theorem}\label{L=4}
If there exists a nontrivial continuous solution of the fractional
 boundary value problem (\ref{L-4}),  then
 \begin{equation}\label{L4} \, \, (i)
 \int_a^b |q(s)|ds> \Big(1+\frac{p}{r}(b-a)\Big)\frac{\Gamma(\alpha)}{(b-a)^{\alpha-1}} \,\, ~~ \text{when} ~~\, \, p>0, ~ \frac{r}{p}>\frac{b-a}{\alpha-1};
   \end{equation}
 {\small  \begin{equation}\label{L4a}  \, \, (ii)
 \int_a^b |q(s)|ds>  \frac{\Gamma(\alpha)}{{\rm \max}\{{\mathcal A}(\alpha, r/p), {\mathcal B}(\alpha, r/p)\}}  \,\, ~~ \text{when} ~~\, \,  p>0, ~ 0\le \frac{r}{p}\le\frac{b-a}{\alpha-1}. \end{equation}}
  \end{theorem}
Using the above Lyapunov-type inequalities, we can find intervals, where linear
combinations of some Mittag-Leffler functions have no real zeros.
 \begin{corollary}
Let $\displaystyle 1<\alpha < 2,~ p>0, ~\frac{r}{p}>\frac{1}{\alpha-1}.$ Then the linear combination of Mittag-Leffler functions given by
$$pE_{\alpha,2}(z)+qrE_{\alpha,1}(z)$$
has no real zeros for
$$z\in \Big(-\Big(1+\frac{p}{r}\Big)\Gamma(\alpha),0\Big].$$
\end{corollary}
This corollary can be established  by considering  the following fractional Sturm-Liouville
eigenvalue problem:
$$
 \left\{\begin{array}{ll}
 ^CD^{\alpha} y(t)+\lambda y(t)=0,  \quad 0<t<1,   \\[0.2cm]
\displaystyle py(0)-ry'(0)=y(1)=0.
 \end{array}
 \right.
 $$
We can apply the foregoing  Lyapunov-type inequalities to study the
nonexistence of solutions for certain fractional boundary value problems. For example, the problem (\ref{L-4}) with $p=1, r=2, a=0, b=1, 3/2 < \alpha <2,$  has no nontrivial solution if
$\int_0^1|q(s)|ds<\frac{3}{2}\Gamma(\alpha).$ As a second example,  there is no nontrivial solution for the problem (\ref{L-4}) with $p=2, r=1, a=0, b=1, 1 < \alpha <2,$ provided that
$$\int_0^1|q(s)|ds<\frac{\Gamma(\alpha)}{\max\{{\mathcal A}(\alpha,1/2), {\mathcal B}(\alpha,1/2)\}}.$$

 In 2015, O'Regan and Samet \cite{ORSa}   obtained  a Lyapunov-type inequality for the   fractional
boundary value problem:
  \begin{equation}\label{L-7}
 \left\{\begin{array}{ll}
D^{\alpha}y(t)+q(t)y(t)=0, \quad a<t<b,\\[0.2cm]
\displaystyle y(a)=y'(a)=y''(a)=y''(b)=0,
 \end{array}
 \right.
 \end{equation}
 where $D^{\alpha}$ is  the standard Riemann-Liouville fractional derivative of fractional order $3<\alpha\le 4$  and  $q: [a,b]\to {\mathbb R}$ is a continuous function.

The integral equation associated with  problem (\ref{L-7}) is
\begin{equation}\label{G-L-7}
y(t)=\int_a^bG(t,s)q(s)y(s)ds,
\end{equation}
where $G(t,s)$ is the Green's function defined by
\begin{equation}\label{Gr-L-7}
G(t,s)=\frac{1}{\Gamma(\alpha)}
\begin{cases}
\displaystyle\frac{(t-a)^{\alpha-1}(b-s)^{\alpha-3}}{(b-a)^{\alpha-3}}, &a\leq t\leq s\leq b,\\[0.2cm]
\displaystyle\frac{(t-a)^{\alpha-1}(b-s)^{\alpha-3}}{(b-a)^{\alpha-3}}-(t-s)^{\alpha-1}, &a\leq s \leq t \leq b.
\end{cases}
\end{equation}
The  Green's function defined in \eqref{Gr-L-7}   satisfies  the following inequality:
$$0\le G(t,s)\le G(b,s)=\frac{(b-s)^{\alpha-3}(s-a)(2b-a-s)}{\Gamma(\alpha)}, ~~ (t,s)\in [a,b]\times [a,b].$$

 The Lyapunov inequality for the problem (\ref{L-7}) is as follows.
 \begin{theorem}\label{L=7}
 If there exists a nontrivial continuous solution of the fractional
 boundary value problem (\ref{L-7}), then
 \begin{equation}\label{L7}
 \int_a^b (b-s)^{\alpha-3}(s-a)(2b-a-s)|q(s)|ds\ge \Gamma(\alpha).
 \end{equation}
 \end{theorem}
 \vskip 0.3cm
 To demonstrate  an application of the above inequality, we consider the eigenvalue problem:
\begin{equation}\label{EIG-1}
 \left\{\begin{array}{ll}
D^{\alpha}y(t)+\lambda y(t)=0,
 \quad 0<t<1,~~ 3<\alpha\le 4,\\
 y(0)=y'(0)=y''(0)=y''(1)=0.
 \end{array}
 \right.
 \end{equation}
\begin{corollary}
If $\lambda$ is an eigenvalue of the problem (\ref{EIG-1}), then
$$|\lambda|\ge \frac{\Gamma(\alpha)}{2B(2,\alpha-2)},$$
where $B$ is the beta function defined by $\displaystyle B(x,y)=\int_0^1 s^{x-1}(1-s)^{y-1}ds, ~ x,y>0.$
\end{corollary}

\vskip 0.3cm
Sitho {\em et al.} \cite{Sit}  established some Lyapunov-type inequalities for the following {\em  hybrid fractional boundary value problem}
\begin{equation}\label{problem-Hybrid}
\begin{cases}
\displaystyle D_a^{\alpha}\left[\displaystyle\frac{y(t)}{f(t,y(t))}-\displaystyle\sum_{i=1}^n I_a^{\beta}h_i(t,y(t))\right]+g(t)y(t)=0, \quad t\in(a,b),\\[4mm]
y(a)=y'(a)=y(b)=0,
\end{cases}
\end{equation}
where $D_a^{\alpha}$ denotes the Riemann-Liouville fractional derivative of order  $\alpha\in(2,3]$ starting from a point $a$, the functions  $g\in L^1((a,b], \mathbb{R})$, $f\in C^1([a,b]\times\mathbb{R}, \mathbb{R}\setminus\{0\})$, $h_i\in C([a,b]\times\mathbb{R}, \mathbb{R})$, $\forall i=1,2,\ldots,n$ and  $I_a^{\beta}$ is the Riemann-Liouville fractional integral of order $\beta\geq \alpha$ with the lower limit at the  point $a$.
We consider two cases:  (I) $h_i=0, i=1,2,\ldots,n$ and (II) $h_i\ne 0, i=1,2,\ldots,n.$

{\bf Case I:} $h_i=0, i=1,2,\ldots,n.$
We consider the problem \eqref{problem-Hybrid} with $h_i(t,\cdot)=0$ for all $t\in [a,b]$. For $\alpha\in(2,3]$, we first construct the Green's function for the following boundary value problem
\begin{equation}\label{problem}
\begin{cases}
D_a^\alpha\left[\displaystyle\frac{y(t)}{f(t,y(t))}\right]+g(t)y(t)=0,\quad t\in(a,b), \\[4mm]
y(a)=y'(a)=y(b)=0,
\end{cases}
\end{equation}
with the assumption that $f$ is continuously differentiable and $f(t,y(t))\neq 0$ for all $t\in[a,b]$.

Let $y\in AC([a,b],\mathbb{R})$ be a solution of the problem \eqref{problem}. Then the function $y$   satisfies the following integral equation
\begin{equation}\label{greensol}
y=f(t,y)\int_a^bG(t,s)g(s)y(s)ds,
\end{equation}
where $G(t,s)$ is the Green's function defined by  (\ref{Gr-L-0}) and satisfies the following properties:
\begin{itemize}
\item[(i)] $G(t,s)\ge 0$, $\forall t,s\in [a,b].$
\item[(ii)]~  $\displaystyle G(t,s)\le H(s):= \frac{(b-s)^{\alpha-1}}{\Gamma(\alpha-1)}.$
\item[(iii)]~ $\displaystyle  \max_{s\in[a,b]} H(s) =\frac{(b-a)^{\alpha-1}}{\Gamma(\alpha-1)}.$
\end{itemize}

\begin{theorem}\label{thmlya}
The necessary condition for existence of a nontrivial solution for the boundary value problem \eqref{problem} is
\begin{equation}\label{lyatype-1}
\frac{\Gamma(\alpha-1)}{\|f\|} \leq \int_a^b (b-s)^{\alpha-1}|g(s)|ds,
\end{equation}
where $\|f\|=\displaystyle\sup_{t\in[a,b], y\in {\mathbb R}}   |f(t,y)| $.
\end{theorem}

{\bf Case II:} $h_i\ne 0, i=1,2,\ldots, n.$

Let $y\in AC[a,b]$ be a solution of the problem \eqref{problem-Hybrid}. Then the function $y$ can be written as
\begin{equation}
y(t)=f(t,y(t))\left[\int_a^b G(t,s)g(s)y(s)ds-\sum_{i=1}^n\int_a^bG^*(t,s)h_i(s,y(s))ds\right],
\end{equation}
where   $G(t,s)$  is defined as in (\ref{Gr-L-0}) and $G^*(t,s)$ is defined by
{\small\begin{equation}\label{greeni}
G^*(t,s)=\frac{1}{\Gamma(\beta)}\begin{cases}
\displaystyle\frac{(b-s)^{\beta-1}(t-a)^{\alpha-1}}{(b-a)^{\alpha-1}}-(t-s)^{\beta-1}, &a\leq s\leq t \leq b,\\[0.3cm]
\displaystyle\frac{(b-s)^{\beta-1}(t-a)^{\alpha-1}}{(b-a)^{\alpha-1}}, &a\leq t\leq s\leq b.
\end{cases}
\end{equation}}

The Green's  function $G^*(t,s),$ which is given by \eqref{greeni}, satisfies the following inequalities:
\begin{itemize}
\item[(i)] $G^*(t,s)\geq0, ~~ \forall t,s\in [a,b];$
\item[(ii)]~ $\displaystyle G^*(t,s)\le J(s):= \frac{(\alpha-1)(b-s)^{\beta-1}}{\Gamma(\beta)}.$
\end{itemize}
 Also we have
 \begin{itemize}
 \item[(iii)]~ $\displaystyle  \max_{s\in[a,b]} J(s) =\frac{(\alpha-1)(b-a)^{\beta-1}}{\Gamma(\beta)}.$
 \end{itemize}

\begin{theorem}\label{thmp201}
Assume that  $|h_i(t,y(t))|\leq |x_i(t)||y(t)|,$ where $x_i\in C([a,b],\mathbb{R})$, $i=1,2,\ldots,n$ and  $[a,b]=[0,1]$. The necessary condition for existence of a nontrivial solution for the problem \eqref{problem-Hybrid} on $[0,1]$ is
\begin{equation}\label{Hlyapunov01}
\Gamma(\alpha-1)\left(\frac{1}{\|f\|}-\frac{(\alpha-1)}{\Gamma(\beta+1)}\sum_{i=1}^n\|x_i\|\right)\leq \int_0^1(1-s)^{\alpha-1}|g(s)|ds.
\end{equation}
\end{theorem}

 \vskip 0.3cm
 In 2016, Al-Qurashi and Ragoub \cite{AlQu}   obtained a Lyapunov-type inequality for a boundary value problem with {\em   integral boundary condition}
 \begin{equation}\label{L-2}
 \left\{\begin{array}{ll}
^CD^{\alpha}y(t)+q(t)y(t)=0, \quad a<t<b,\\[0.2cm]
\displaystyle y(a)+\mu\int_a^by(s)q(s)ds=y(b),
 \end{array}
 \right.
 \end{equation}
 where $^CD^{\alpha}$ is the Caputo fractional derivative of order  $0<\alpha\le 1,$  $q: [a,b]\to {\mathbb R}$ is a continuous function, $a, b$ are consecutive zeros of the solution $y$ and $\mu$  is positive.

The function $y\in C([a,b], {\mathbb R})$ is a solution of the boundary value problem (\ref{L-2})
if  and only if  $y$ satisfies the integral equation
\begin{equation}\label{G-L-2}
y(t)=\int_a^bG(t,s)q(s)y(s)ds,
\end{equation}
where $G(t,s)$ is the Green's function defined by
\begin{equation}\label{Gr-L-2}
G(t,s)=\frac{1}{\Gamma(\alpha)}
\begin{cases}
\displaystyle \frac{(b-s)^{\alpha-1}}{\mu(b-a)}-\frac{(b-s)^{\alpha}}{(b-a)\alpha}, &a\leq t\leq s\leq b,\\[0.3cm]
\displaystyle \frac{(b-s)^{\alpha-1}}{\mu(b-a)}-\frac{(b-s)^{\alpha}}{(b-a)\alpha}-(t-s)^{\alpha-1}, &a\leq s \leq t \leq b,
\end{cases}
\end{equation}
and satisfies the following properties:
 \begin{itemize}
 \item[(i)] $G(t,s)\ge 0, ~~\mbox{for all}~~ a\le t,s\le b;$
 \item[(ii)]~ $\max_{t\in [a,b]}G(t,s)=G(b,s), ~ s\in [a,b];$
 \item[(iii)]~ $G(b,s)$ has a unique maximum given  by {\small  $\displaystyle  \max_{s\in [a,b]}G(b,s)=\frac{\alpha(b-a+\mu)(\alpha\mu+1)^{\alpha-1}}{\Gamma(\alpha)\mu^{\alpha}(b-a)},$}
 provident $0<\mu(b-a)<\alpha.$
 \end{itemize}
 We describe the Lyapunov's inequality for the problem (\ref{L-2}) as follows.
 \begin{theorem}\label{L=2}
The
 boundary value problem (\ref{L-2})  has a nontrivial solution
provided that the real and continuous function $q$ satisfies the
following inequality
 \begin{equation}\label{L2}
 \int_a^b |q(s)|ds>\frac{\Gamma(\alpha)\mu^{\alpha}(b-a)}{\alpha(b-a+\mu)(\alpha\mu+1)^{\alpha-1}}.
 \end{equation}
 \end{theorem}

 In 2016, Fereira \cite{Fer}   obtained a  Lyapunov-type inequality for a {\em sequential fractional
boundary value problem}
  \begin{equation}\label{L-5}
 \left\{\begin{array}{ll}
(D^{\alpha}D^{\beta}y)(t)+q(t)y(t)=0, \quad a<t<b,\\[0.2cm]
\displaystyle y(a)=y(b)=0,
 \end{array}
 \right.
 \end{equation}
 where $D^{\delta},~\delta=\alpha, \beta$ stands for the Riemann-
Liouville  fractional derivative and  $q: [a,b]\to {\mathbb R}$ is a continuous function.
Assuming that
(\ref{L-5})  has a nontrivial solution $y\in C[a, b]$ of the form
$$y(t)=c\frac{\Gamma(\alpha)}{\Gamma(\alpha+\beta)}(t-a)^{\alpha+\beta-1}-(I^{\alpha+\beta}qy)(t),$$
it follows by Proposition 2.6 and the fact $\displaystyle I^{\beta}(t-a)^{\alpha-1}= \frac{\Gamma(\alpha)}{\Gamma(\alpha+\beta)}(t-a)^{\alpha+\beta-1}$ that $y'$ is integrable in $[a, b].$  Then,
as argued in (\cite[Section 2.3.6-2.3.7]{Kil}), we have
$$(D^{\alpha}D^{\beta}y)(t)=(D^{\alpha+\beta}y)(t).$$
 The following result is therefore an immediate consequence of Theorem \ref{L=0}.
 \begin{theorem}\label{L=5} (Riemmann-Liouville Case)
Let $0<\alpha, \beta\le 1$ with $1<\alpha+\beta\le 2.$ If there exists a nontrivial continuous solution of the fractional  boundary value problem (\ref{L-5}), then
 \begin{equation}\label{L5}
 \int_a^b |q(s)|ds>\Gamma(\alpha+\beta)\Big(\frac{4}{b-a}\Big)^{\alpha+\beta-1}.
 \end{equation}
 \end{theorem}

 As an application we consider the following sequential fractional differential equation
 \begin{equation}\label{exam-2}
 (D^{\alpha}D^{\alpha}y)(t)+\lambda^2 y(t)=0, ~~ \lambda\in {\mathbb R}, ~~ t\in (0,1), ~~ \frac{1}{2}<\alpha\le 1.
 \end{equation}
 The fundamental set of solutions to (\ref{exam-2}) is
 $$\{\cos_{\alpha}(\lambda t), \sin_{\alpha}(\lambda t)\},$$
 where
  {\small $$\cos_{\alpha}(\lambda t)=\sum_{j=0}^{\infty}(-1)^j\lambda^{2j}\frac{t^{(2j+1)\alpha-1}}{\Gamma((2j+1)\alpha)} ~~ \mbox{and}~~ \sin_{\alpha}(\lambda t)=\sum_{j=0}^{\infty}(-1)^j\lambda^{2j+1}\frac{t^{(j+1)2\alpha-1}}{\Gamma((j+1)2\alpha)}.$$}
 Therefore the general solution of (\ref{exam-2}) can be written as
$$y(t) =c_1\cos_{\alpha}(\lambda t)+c_2\sin_{\alpha}(\lambda t), ~~~ c_1, c_2\in {\mathbb R}.$$
Now, the nontrivial solutions of (\ref{exam-2}) for which the boundary conditions $y(0) = 0 =
y(1)$ hold, satisfy $\sin_{\alpha}(\lambda)=0,$ where $\lambda$ is a real number different from zero (the eigenvalue of the problem). By Theorem \ref{L=5}, the following
inequality then holds
 $$\lambda^2>\Gamma(2\alpha)4^{2\alpha-1},$$
which can alternatively be expressed in form of the following result.
\begin{corollary}
 Let $\displaystyle \frac{1}{2}<\alpha\le 1.$ If
$$|t|\le \sqrt{\Gamma(2\alpha)4^{2\alpha-1}}, ~~ t\ne 0,$$
then $\sin_{\alpha}(t)$ has no real zeros.
\end{corollary}

 In \cite{Fer},  Fereira replaced the  Riemmann-Liouville fractional derivative in the problem (\ref{L-5}) with the  Caputo  fractional derivative and obtained the following Lyapunov-type inequality:
 \begin{equation}\label{L5-a}
 \int_a^b |q(s)|ds>\frac{\Gamma(\alpha+\beta)}{(b-a)^{\alpha+\beta-1}}\frac{(\alpha+2\beta-1)^{\alpha+2\beta-1}}{(\alpha+\beta-1)^{\alpha+\beta-1}\beta^{\beta}}.
 \end{equation}

In 2016, Al-Qurashi and Ragoub \cite{AlQu-1}   obtained a Lyapunov-type inequality for a {\em  fractional
boundary value problem: }
  \begin{equation}\label{L-9}
 \left\{\begin{array}{ll}
D^{\alpha}y(t)+q(t)y(t)=0, \quad 1<t<e,\\[0.2cm]
\displaystyle y(a)=y(b)=y''(a)=y''(b)=0,
 \end{array}
 \right.
 \end{equation}
 where $D^{\alpha}$ is  the standard Riemann-Liouville fractional derivative of   order $3<\alpha\le 4$  and  $q: [a,b]\to {\mathbb R}$ is a continuous function.

 The function $y\in C([a,b], {\mathbb R})$ is a solution of the boundary value problem (\ref{L-9})
if and only if $y$ satisfies the integral equation
\begin{equation}\label{G-L-9}
y(t)=\int_a^bG(t,s)q(s)y(s)ds,
\end{equation}
where $G(t,s)$ is the Green's function defined by
{\small\begin{equation}\label{Gr-L-9}
G(t,s)=\frac{1}{\Gamma(\alpha)}
\begin{cases}
\displaystyle -(t-a)(b-s)^{\alpha-1}\\
\displaystyle +\frac{\alpha(\alpha-1)(b-a)(b-s)^{\alpha-3}}{6}(t-a)\Big[1-\frac{(t-a)^2}{(b-a)^2}\Big], &a\leq t\leq s\leq b,\\[0.4cm]
\displaystyle(t-s)^{\alpha-1}-(t-a)(b-s)^{\alpha-1}\\
\displaystyle+\frac{\alpha(\alpha-1)(b-a)(b-s)^{\alpha-3}}{6}(t-a)\Big[1-\frac{(t-a)^2}{(b-a)^2}\Big], &a\leq s \leq t \leq b,
\end{cases}
\end{equation}}
and  satisfies the relation:
$$0\le G(t,s)\le G(b,s)=\frac{(1-(b-a))(b-a)^{\alpha-1}}{\Gamma(\alpha)},~~ (t,s)\in [a,b]\times [a,b].$$

The Lyapunov-type inequality for the problem (\ref{L-9}) is given in the following result.

 \begin{theorem}\label{L=9}
 If there exists a nontrivial continuous solution to the fractional
 boundary value problem (\ref{L-9}), then
 \begin{equation}\label{L9}
 \int_a^b  |q(s)|ds>\frac{\Gamma(\alpha)}{(1-(b-a))(b-a)^{\alpha-1}}.
 \end{equation}
 \end{theorem}

 In order to illustrate Theorem \ref{L=9}, we apply the Lyapunov-type inequality (\ref{L9}) to find a bound for $\lambda$ so that the following eigenvalue problem has a nontrivial solution:
\begin{equation}\label{exam-3}
 \left\{\begin{array}{ll}
D^{\alpha}y(t)+\lambda y(t)=0,  \quad 0<t<\frac{1}{2},~ 3<\alpha\le 4,\\[0.2cm]
\displaystyle y(1)=y\Big(\frac{1}{2}\Big)=y''(1)=y''\Big(\frac{1}{2}\Big)=0.
 \end{array}
 \right.
 \end{equation}
  \begin{corollary}
 If $\lambda$ is an eigenvalue of the  fractional boundary value problem (\ref{exam-3}), then the following inequality
holds
 $$|\lambda|\ge \frac{\Gamma(\alpha)}{2^{-\alpha}}.$$
 \end{corollary}

 \vskip 0.3cm

In 2016, Dhar {\em et al.} \cite{DKM} derived Lyapunov-type inequalities  for the following boundary value problem
 \begin{equation}\label{exam-15}
 \left\{\begin{array}{ll}
 D^{\alpha}y(t)+q(t) y(t)=0,  \quad  a<t<b, ~ 1<\alpha\le 2, \\[0.2cm]
\displaystyle D^{\alpha-2}y(a)=D^{\alpha-2}y(b)=0,
 \end{array}
 \right.
 \end{equation}
where $D^{\alpha}$ is the Riemann-Liouville fractional derivative of order $\alpha$ ($1<\alpha\le 2$), $q\in L( [a,b], {\mathbb R}).$
Their main result on fractional Lyapunov-type inequalities is the following.
 \begin{theorem}
 \begin{itemize}
\item[(a)]~ If the problem (\ref{exam-15}) has a nontrivial solution, then
$$\max_{t\in [a,b]}\Bigg\{\int_a^b|D^{2-\alpha}[G(t,s)q(s)]|ds\Bigg\}>1,$$
where $D^{2-\alpha}[G(t,s)q(s)]$ is  the right-sided fractional derivative of $G(t, s)q(s)$ with respect to
$s$ with
  $$G(t,s)=\frac{1}{b-a }
\begin{cases}
\displaystyle (t-a)(b-s), &a\leq t\leq s\leq b,\\[0.2cm]
\displaystyle (s-a)(b-t), &a\leq s \leq t \leq b.
\end{cases}
$$
\item[(b)]~  If  problem (\ref{exam-15}) has a nontrivial solution and $D^{\alpha-2}y(t)\ne 0$ on $(a,b),$ then
$$\max_{t\in [a,b]}\Bigg\{\int_a^b\Big[D^{2-\alpha}[G(t,s)q(s)]\Big]_{+}ds\Bigg\}>1,$$
where $D^{2-\alpha}[G(t,s)q(s)]_{+}$ is  the the positive part of $D^{2-\alpha}[G(t,s)q(s)].$
\end{itemize}
\end{theorem}
As a special case we have the following corollary.
\begin{corollary}
Assume that $D^{2-\alpha}_{b-}[G(t,s)q(s)]\ge 0$ for $t,s\in [a,b]$ so that the problem (\ref{exam-15}) has a nontrivial solution. Then
$$\int_a^bq_+(t)dt>\frac{\alpha^{\alpha}\Gamma(\alpha-1)}{(\alpha-1)^{\alpha-1}(b-a)^{\alpha-1}}.$$
\end{corollary}
Next  we  consider the sequential fractional boundary value problem
 \begin{equation}\label{exam-15a}
 \left\{\begin{array}{ll}
 \Big[\big(D^{\beta}_{a+}(D^{\alpha} y)\big)\Big](t)+q(t) y(t)=0,  \quad a<t<b, ~ 0<\alpha,\beta\le 1, \\[0.2cm]
\displaystyle (D^{\alpha-1}y)(a^+)=(D^{\alpha-1}y)(b)=0,
 \end{array}
 \right.
 \end{equation}
 which is equivalent to the integral equation
 $$y(t)=\int_a^bG(t,s)q(s)D^{1-\alpha}_{a+}y(s)ds,$$
 where
 $$G(t,s)=\frac{1}{\Gamma(\beta+1)}
\begin{cases}
\displaystyle \frac{(t-a)^{\beta}(b-s)^{\beta}}{(b-a)^{\beta}}, &a\leq t\leq s\leq b,\\[0.4cm]
\displaystyle \frac{(t-a)^{\beta}(b-s)^{\beta}}{(b-a)^{\beta}}-(t-s)^{\beta}, &a\leq s \leq t \leq b.
\end{cases}
$$
In the following result, we express the fractional Lyapunov-type inequalities for  problem (\ref{exam-15a}).
 \begin{theorem}
 \begin{itemize}
\item[(a)]~ If  problem (\ref{exam-15a}) has a nontrivial solution, then
$$\max_{t\in [a,b]}\Bigg\{\int_a^b|D^{1-\alpha}[G(t,s)q(s)]|ds\Bigg\}>1.$$
\item[(b)]~  If  problem (\ref{exam-15a}) has a nontrivial solution and  $(D^{\alpha-1}_{a+}y)(t)\ne 0$ on $(a,b),$ then
$$\max_{t\in [a,b]}\Bigg\{\int_a^b\Big[D^{1-\alpha}[G(t,s)q(s)]\Big]_{+}ds\Bigg\}>1.$$
\end{itemize}
\end{theorem}
As a special case we have  the following corollary.
\begin{corollary}
Assume that $D^{1-\alpha}[G(t,s)q(s)]\ge 0$ for $t,s\in [a,b],$ $1<\alpha+\beta\le 2$ and  the problem (\ref{exam-15a}) has a nontrivial solution. Then
$$\int_a^bq_+(t)dt>\frac{(\alpha+2\beta-1)^{\alpha+2\beta-1}\Gamma(\alpha)\Gamma(\beta+1)}{(\alpha+\beta-1)^{\alpha+\beta-1}\beta^{\beta}(b-a)^{\alpha+\beta-1}}.$$
\end{corollary}

For some similar results on  fractional boundary value problems of order $\alpha \in (2, 3]$, see \cite{DK-EJDE}.
\vskip 0.3cm

 In 2017, Jleli {\em et al.} \cite{L9}   obtained   Lyapunov-type inequality for {\em higher order fractional
boundary value problem }
  \begin{equation}\label{L-8}
 \left\{\begin{array}{ll}
D^{\alpha}y(t)+q(t)y(t)=0, \quad a<t<b,\\[0.2cm]
\displaystyle y(a)=y'(a)=\ldots =y^{(n-2)}(a)=0,~ y(b)=I^{\alpha}(yh)(b),
 \end{array}
 \right.
 \end{equation}
 where $n\in {\mathbb N}, n-1<\alpha<n,$ $D^{\alpha}$ is  the standard Riemann-Liouville fractional derivative of  order $\alpha,$
 $I^{\alpha}$ denotes the Riemann-Liouville fractional integral of order $\alpha,$ and  $q, h: [a,b]\to {\mathbb R}$ are continuous functions.

 The function  $y$ is a solution of the boundary value problem (\ref{L-8})
if  and only if  $y$ satisfies the integral equation
\begin{equation}\label{G-L-8}
y(t)=\int_a^bG(t,s)(q(s)+h(s))y(s)ds+\frac{1}{\Gamma(\alpha)}\int_a^t(t-s)^{\alpha-1}h(s)y(s)ds,
\end{equation}
where $G(t,s)$ is the Green's function given by (\ref{Gr-L-0})
 such that
$$0\le G(t,s)\le G(s^*,s)=\frac{(s-a)^{\alpha-1}(b-s)^{\alpha-1}}{\Gamma(\alpha)(b-a)^{\alpha-1}\Big[1-\Big(\frac{b-s}{b-a}\Big)^{\frac{\alpha-1}{\alpha-2}}\Big]^{\alpha-2}}, ~a<s<b,$$
with
$$s^*=\frac{s-a\Big(\frac{b-s}{b-a}\Big)^{\frac{\alpha-1}{\alpha-2}}}{1-\Big(\frac{b-s}{b-a}\Big)^{\frac{\alpha-1}{\alpha-2}}}.$$
 The following result presents the  Lyapunov-type inequality for  problem (\ref{L-8}).
 \begin{theorem}\label{L=8}
 Let $n\in {\mathbb N}$ with $n\ge 3.$  If $y$
 is a nontrivial   solution of the fractional
 boundary value problem (\ref{L-8}), then
 {\small  \begin{equation}\label{L8}
 \int_a^b\Bigg(|q(s)+h(s)|+\frac{\Big(1-z_{\alpha}^{\frac{\alpha-1}{\alpha-2}}\Big)^{\alpha-2}}{z_{\alpha}^{\alpha-1}(1-z_{\alpha})^{\alpha-1}}|h(s)|\Bigg)ds\ge \frac{\Gamma(\alpha)}{(b-a)^{\alpha-1}}\frac{\Big(1-z_{\alpha}^{\frac{\alpha-1}{\alpha-2}}\Big)^{\alpha-2}}{z_{\alpha}^{\alpha-1}(1-z_{\alpha})^{\alpha-1}},
 \end{equation}}
 where $z_{\alpha}$ is the unique zero of the nonlinear algebraic equation
 $$z^{\frac{2\alpha-3}{\alpha-2}}-2z+1=0$$
 in the interval $\Big(0,\Big(\frac{2\alpha-4}{2\alpha-3}\Big)^{\frac{\alpha-2}{\alpha-1}}\Big).$
 \end{theorem}

 \begin{theorem}\label{L=8a}
 Let $n=2.$   If $y$
 is a nontrivial   solution of the fractional
 boundary value problem
  \begin{equation}\label{L-8-2}
 \left\{\begin{array}{ll}
D^{\alpha}y(t)+q(t)y(t)=0, \quad a<t<b,\\[0.2cm]
\displaystyle y(a)=0,~ y(b)=I^{\alpha}(hy)(b),
 \end{array}
 \right.
 \end{equation}
  then
 \begin{equation}\label{L8a}
 \int_a^b\Bigg(|q(s)+h(s)|+4^{\alpha-1}h(s)|\Bigg)ds\ge \Gamma(\alpha)\Big(\frac{4}{b-a}\Big)^{\alpha-1}.
 \end{equation}
 \end{theorem}

 \vskip 0.3cm

In 2017, Cabrera {\em et al.}  \cite{CLS} obtained  a Lyapunov-type inequality for a sequential fractional
boundary value problem
  \begin{equation}\label{L-11}
 \left\{\begin{array}{ll}
^CD^{\alpha}y(t)+q(t)y(t)=0, \quad a<t<b,~~ \alpha\in (n-1,n], ~ n\ge 4,\\[0.2cm]
\displaystyle y^i(a)=y''(b)=0,~~ 0\le i\le n-1, ~ i\ne 2,
 \end{array}
 \right.
 \end{equation}
 where $^CD^{\alpha}$ is  the Caputo fractional derivative of fractional order $\alpha\ge 0$  and  $q: [a,b]\to {\mathbb R}$ is a continuous function.

 The function $y\in C([a,b], {\mathbb R})$ is a solution of the boundary value problem (\ref{L-11})
if and only if it  satisfies the integral equation
\begin{equation}\label{G-L-11}
y(t)=\int_a^bG(t,s)q(s)y(s)ds,
\end{equation}
where $G(t,s)$ is the Green's function defined by
 {\small \begin{equation}\label{Gr-L-11}
G(t,s)=\frac{1}{\Gamma(\alpha)}
\begin{cases}
\displaystyle \frac{1}{2}(\alpha-1)(\alpha-2)(t-a)^2(b-s)^{\alpha-3}-(t-s)^{\alpha-1}, &a\leq s\leq t\leq b,\\[0.3cm]
\displaystyle \frac{1}{2}(\alpha-1)(\alpha-2)(t-a)^2(b-s)^{\alpha-3}, &a\le t\le  s  \leq b,
\end{cases}
\end{equation}}
such that  $G(t,s)\ge 0$ for $t,s\in [a,b],$ $|G(t,s)|\le G(b,s)$ for $t,s\in [a,b]$ and
{\small $$|G(t,s)|\le G(b,s)=\frac{1}{2}(\alpha-1)(\alpha-2)(b-a)^2(b-s)^{\alpha-3}-(b-s)^{\alpha-1}, ~~ (t,s)\in [a,b]\times [a,b].$$}

 Their result is  as follows.
 \begin{theorem}\label{L=11}
 If there exists a nontrivial continuous solution of the fractional
 boundary value problem (\ref{L-11}), then
 \begin{equation}\label{L11}
 \int_a^b \Big[\frac{1}{2}(\alpha-1)(\alpha-2)(b-a)^2(b-s)^{\alpha-3}-(b-s)^{\alpha-1}\Big] |q(s)|ds\ge \Gamma(\alpha).
 \end{equation}
 \end{theorem}

\vskip 0.3cm

In 2017,  Wang {\em et al.} \cite{WPL}  obtained a  Lyapunov-type inequality for the {\em higher order fractional
boundary value problem }
  \begin{equation}\label{Oct-1}
 \left\{\begin{array}{ll}
D^{\alpha}y(t)+q(t)y(t)=0, \quad a<t<b,\\[0.2cm]
\displaystyle y(a)=y'(a)=\ldots =y^{(n-2)}(a)=0,~ y^{(n-2)}(b)=0,
 \end{array}
 \right.
 \end{equation}
 where $n\in {\mathbb N}, 2<n-1<\alpha\le n,$ $D^{\alpha}$ is  the standard Riemann-Liouville fractional derivative of  order $\alpha,$
  and  $q: [a,b]\to {\mathbb R}$ is a continuous function.

 The function  $y$ is a solution of the boundary value problem (\ref{Oct-1})
if  and only if  $y$ satisfies the integral equation
\begin{equation*}
y(t)=\int_a^bG(t,s)q(s)y(s)ds,
\end{equation*}
where $G(t,s)$ is the Green's function given by
$$G(t,s)=\frac{1}{\Gamma(\alpha)}
\begin{cases}
\displaystyle \frac{(t-a)^{\alpha-1}(b-s)^{\alpha-n+1}}{(b-a)^{\alpha-n+1}}, &a\leq t\leq s\leq b,\\[0.4cm]
\displaystyle\frac{(t-a)^{\alpha-1}(b-s)^{\alpha-n+1}}{(b-a)^{\alpha-n+1}}-(t-s)^{\alpha-1}, &a\leq s \leq t \leq b,
\end{cases}
$$
 such that
$$0\le G(t,s)\le G(b,s)=\frac{(b-s)^{\alpha-n+1}(s-a)}{\Gamma(\alpha)}\sum_{i=1}^{n-2}(-1)^{i-1}C_{n-2}^i(b-a)^{n-2-i}(s-a)^{i-1},$$
$(t,s)\in [a,b]\times [a,b]$ and  $C_{n-2}^i$ is the binomial coefficient.

Their Lyapunov-type inequality for the problem  (\ref{Oct-1}) is given in the following theorem.
\begin{theorem}\label{t-Oct-1}
If there exists a nontrivial continuous solutions $y$ of the fractional boundary value problem (\ref{Oct-1}), and $q$ is a real continuous function, then
$$\int_a^b(b-s)^{\alpha-n+1}(s-a)\sum_{i=1}^{n-2}(-1)^{i-1}C_{n-2}^i(b-a)^{n-2-i}(s-a)^{i-1}|q(s)|ds\ge \Gamma(\alpha).$$
\end{theorem}

\begin{corollary}
If the fractional boundary value problem  (\ref{Oct-1})  has a nontrivial continuous solution, then
$$\int_a^b|q(s)|ds\ge\frac{\Gamma(\alpha)(\alpha-n+2)^{\alpha-n+2}}{(n-2)(\alpha-n+1)^{\alpha-n+1}(b-a)^{\alpha-1}}.$$
\end{corollary}

The following result shows the application of the above Lyapunov-type inequality to eigenvalue problems.

\begin{corollary}
If $\lambda$  is an eigenvalue to the fractional boundary value problem
\begin{equation*}
 \left\{\begin{array}{ll}
D^{\alpha}y(t)+\lambda y(t)=0, \quad a<t<b,\\[0.2cm]
\displaystyle y(a)=y'(a)=\ldots =y^{(n-2)}(a)=0,~ y^{(n-2)}(b)=0,
 \end{array}
 \right.
 \end{equation*}
 then
 $$|\lambda|\ge \frac{\Gamma(\alpha)(\alpha-n+3)(\alpha-n+2)}{n-2}.$$
\end{corollary}

In 2017, Fereira \cite{r5}  obtained a Lyapunov-type inequality for so-called  {\em anti-periodic}  boundary value problem:
 \begin{equation}\label{L-12}
 \left\{\begin{array}{ll}
^CD^{\alpha}y(t)+q(t)y(t)=0, \quad a<t<b,\\[0.2cm]
 y(a)+y(b)=0=y'(a)+y'(b),
 \end{array}
 \right.
 \end{equation}
 where $^CD^{\alpha}$ is the Caputo fractional derivative of order  $1<\alpha\le 2$ and $q: [a,b]\to {\mathbb R}$ is a continuous function.

Then $y\in C([a,b], {\mathbb R})$ is a solution of the boundary value problem (\ref{L-12})
if and only if it satisfies the integral equation
\begin{equation}\label{G-L-12}
y(t)=\int_a^b (b-s)^{\alpha-2}H(t,s)q(s)y(s)ds,
\end{equation}
where $H(t,s)$ is  defined by
 {\small \begin{equation}\label{Gr-L-12}
H(t,s)=\frac{1}{\Gamma(\alpha)}
\begin{cases}
\displaystyle \Big(-\frac{b-a}{4}+\frac{t-a}{2}\Big)(\alpha-1)+\frac{b-s}{2}, &a\leq t\leq s\leq b,\\[0.3cm]
\displaystyle \Big(-\frac{b-a}{4}+\frac{t-a}{2}\Big)(\alpha-1)+\frac{b-s}{2}-\frac{(t-s)^{\alpha-1}}{(b-a)^{\alpha-2}}, &a\leq s \leq t \leq b.
\end{cases}
\end{equation}}
 Here the  function $H$ satisfies the following property:
 $$|H(t,s)|\le \frac{(b-a)(3-\alpha)}{4}, ~~ (t,s)\in [a,b]\times [a,b].$$

The Lyapunov-type inequality for the problem (\ref{L-12}) is given in the following result.

 \begin{theorem}\label{L=12}
 If (\ref{L-12}) admits a nontrivial continuous solution,
then
 \begin{equation}\label{L12}
 \int_a^b(b-s)^{\alpha-2} |q(s)|ds\ge \frac{4}{(b-a)(3-\alpha)}.
 \end{equation}
 \end{theorem}
 Inequality (\ref{L12}) is useful in finding  a bound for the
possible eigenvalues of the fractional boundary value problem:
$$
 \left\{\begin{array}{ll}
D^{\alpha}y(t)+\lambda y(t)=0,~ a<t<b,\\
 y(a)+y(b)=0=y'(a)+y'(b),
  \end{array}
 \right.
$$
that is,  an eigenvalue $\lambda\in {\mathbb R}$ satisfies the inequality
$$|\lambda|>\frac{4(\alpha-1)}{(b-a)^{\alpha}(3-\alpha)}.$$

\vskip 0.3cm

In 2017, Agarwal and  Zbekler \cite{AgZb} obtained a Lyapunov-type inequality for the following  fractional boundary value problem with  {\em Dirichlet boundary conditions}
\begin{equation}\label{L-N-1}
 \left\{\begin{array}{ll}
(D^{\alpha}y)(t)+p(t)|y(t)|^{\mu-1}y(t)+q(t)|y(t)|^{\gamma-1}y(t)=f(t), \quad a<t<b,\\[0.2cm]
\displaystyle y(a)=0, ~~ y(b)=0,
 \end{array}
 \right.
 \end{equation}
where  $D^{\alpha}$ is the Riemann-Liouville fractional derivative, $p, q, f\in C[t_0,\infty)$ and $0<\gamma<1<\mu<2.$ No sign restrictions are imposed on the potentials $p$ and $q,$ and the forcing term
$f.$

The problem (\ref{L-N-1}) is equivalent to the following integral equation
$$y(t)=\int_a^bG(t,s)[p(s)y^{\mu}(s)+q(s)y^{\gamma}(s)-f(s)]ds,$$
where $G(t,s)$ is the Green's function defined by (\ref{Gr-L-0}).
\vskip 0.3cm
Their Lyapunov-type inequality for the problem (\ref{L-N-1}) is as follows.
\begin{theorem}\label{t-N-1}
Let $y$ be a nontrivial solution of the problem (\ref{L-N-1}). If $y(t)\ne 0$ in $(a,b),$ then the inequality
 {\small \begin{equation}\label{LN1}
\Big(\int_a^b[p^+(t)+q^+(t)]dt\Big)\Big(\int_a^b[\mu_0 p^+(t)+\gamma_0 q^+(t)+|f(t)|]dt\Big)>\frac{4^{2\alpha-3}\Gamma^2(\alpha)}{(b-a)^{2\alpha-2}}
\end{equation}}
holds, where $u^+=\max\{u,0\}, u=p,q$ and
$$\mu_0=(2-\mu)\mu^{\mu/(2-\mu)}2^{2/(\mu-2)}>0,~~
\gamma_0=(2-\gamma)\gamma^{\gamma/(2-\gamma)}2^{2/(\gamma-2)}>0.$$
\end{theorem}

In 2017,  Zhang and  Zheng \cite{ZZ} considered the Riemann-Liouville fractional differential equations with mixed nonlinearities of order $\alpha\in (n-1,n]$ for $n\ge 3$
\begin{equation}\label{sep-1}
(D^{\alpha}y)(t)+p(t)|y(t)|^{\mu-1}y(t)+q(t)|y(t)|^{\gamma-1}y(t)=f(t),
\end{equation}
where $p,q,f\in C([t_0,\infty), {\mathbb R})$ and the constants satisfy $0<\gamma<1<\mu<n~ (n\ge 3).$ Equation (\ref{sep-1}) subjects to the following two kinds boundary conditions, respectively:
\begin{equation}\label{sep-bc-1}
y(a)=y'(a)=y''(a)=\ldots=y^{(n-2)}(a)=y(b)=0,
\end{equation}
and
\begin{equation}\label{sep-bc-2}
y(a)=y'(a)=y''(a)=\ldots=y^{(n-2)}(a)=y'(b)=0,
\end{equation}
where $a$ and $b$ are two consecutive zeros of the function $y.$

Obviously, it is easy to see that equation (\ref{sep-1}) has two special forms; one is the forced {\em sub-linear} ($p(t)=0$) fractional equation
\begin{equation}\label{sep-1-a}
(D^{\alpha}y)(t)+q(t)|y(t)|^{\gamma-1}y(t)=f(t),~ 0<\gamma<1,
\end{equation}
and the other is the forced {\em super-linear} ($q(t)=0$) fractional equation
\begin{equation}\label{sep-1-b}
(D^{\alpha}y)(t)+q(t)|y(t)|^{\gamma-1}y(t)=f(t),~ 1<\mu<n.
\end{equation}

Their Lyapunov-type inequalities for the problems (\ref{sep-1})-(\ref{sep-bc-1}) and (\ref{sep-1})-(\ref{sep-bc-2}) are respectively the following:
\begin{theorem}\label{sep-t-1}
Let $y$ be a positive solution of the boundary value problem (\ref{sep-1})-(\ref{sep-bc-1}) in $(a,b).$ Then
\begin{eqnarray*}
&&\Bigg(\int_a^b[(b-s)(s-a)]^{\alpha-1}\Bigg[1-\Bigg(\frac{b-s}{b-a}\Bigg)^{\frac{\alpha-1}{\alpha-2}}\Bigg]^{2-\alpha}[\mu_0p^+(s)+\gamma_0q^+(s)+f^{-}(s)]ds\Bigg)\\
&&\times\Bigg(\int_a^b[(b-s)(s-a)]^{\alpha-1}\Bigg[1-\Bigg(\frac{b-s}{b-a}\Bigg)^{\frac{\alpha-1}{\alpha-2}}\Bigg]^{2-\alpha}[p^+(s)+q^+(s)]ds\Bigg)\\
&&>[\Gamma(\alpha)(b-a)^{\alpha-1}]^{\frac{n}{n-1}}(n-1)n^{\frac{n}{1-n}},
\end{eqnarray*}
where $\displaystyle\mu_0=(n-\mu)\mu^{\frac{\mu}{n-\mu}}n^{\frac{n}{n-\mu}}$ and $\displaystyle\gamma_0=(n-\gamma)\gamma^{\frac{\gamma}{n-\gamma}}n^{\frac{n}{\gamma-n}}.$
\end{theorem}

\begin{theorem}\label{sep-t-2}
Let $y$ be a positive solution of the boundary value problem (\ref{sep-1})-(\ref{sep-bc-2}) in $(a,b).$ Then
\begin{eqnarray*}
&&\Bigg(\int_a^b[(b-s)^{\alpha-2}(s-a)][\mu_0p^+(s)+\gamma_0q^+(s)+f^{-}(s)]ds\Bigg)\\
&&\times\Bigg([b-s)^{\alpha-2}(s-a)][p^+(s)+q^+(s)]ds\Big)^{\frac{1}{\mu-1}}\\
&&>\Gamma(\alpha)^{\frac{n}{n-1}}(n-1)n^{\frac{n}{1-n}},
\end{eqnarray*}
where $\mu_0$ and $\gamma_0$ are the same as in Theorem \ref{sep-t-2}.
\end{theorem}

In 2017,  Chidouh and Torres \cite{ChTo} extended
the linear term $q(t)y(t)$ to a nonlinear term of the form $q(t)f (y(t))$ and   obtained a generalized Lyapunov's inequality for the fractional boundary value problem
\begin{equation}\label{L-N-3}
 \left\{\begin{array}{ll}
(D^{\alpha}y)(t)+q(t)f(y(t))=0, \quad a<t<b,~~ 1<\alpha\le 2,\\[0.2cm]
\displaystyle y(a)=y(b)=0,
 \end{array}
 \right.
 \end{equation}
where  $D^{\alpha}$ is the Riemann-Liouville fractional derivative of order $\alpha,$  $f: {\mathbb R}^+\to {\mathbb R}^+$ is continuous and $q: [a,b]\to {\mathbb R}^+$ is a Lebesgue integrable  function.

An integral equation equivalent to the Problem (\ref{L-N-3}) is
$$y(t)=\int_a^bG(t,s)q(s)f(y(s))ds,$$
where $G(t,s)$ is the Green's function defined by (\ref{Gr-L-0}) and satisfies the following properties:
\begin{itemize}
\item[(i)] $G(t,s)\ge 0, $ for all $(t,s)\in [a,b]\times [a,b];$
\item[(ii)]~ $\max_{t\in [a,b]}G(t,s)=G(s,s), ~s\in [a,b];$
\item[(iii)]~ $G(s,s)$ has a unique maximum given by $\displaystyle \max_{s\in [a,b]}G(s,s)=\frac{(b-a)^{\alpha-1}}{4^{\alpha-1}\Gamma(\alpha)}.$
\end{itemize}
The Lyapunov-type inequality for the problem (\ref{L-N-3}) is as follows.
\begin{theorem}\label{t-N-3}
Let $q$ be a real nontrivial Lebesgue integrable  function. Assume that $f\in C({\mathbb R}^+, {\mathbb R}^+)$ is a concave and nondecreasing function. If the fractional boundary value problem (\ref{L-N-3})  has a nontrivial solution $y,$ then
\begin{equation}\label{LN3}
\int_a^b q(s)ds>\frac{4^{\alpha-1}\Gamma(\alpha)\eta}{(b-a)^{\alpha-1}f(\eta)}.
\end{equation}
where $\eta=\max_{t\in [a,b]}y(t).$
\end{theorem}
\vskip 0.3cm
In 2016, Ma \cite{Ma}   obtained a generalized form of  Lyapunov's inequality for the fractional boundary value problem
\begin{equation}\label{L-N-4}
 \left\{\begin{array}{ll}
(D^{\alpha}y)(t)+q(t)f(y(t))=0, \quad a<t<b,~~ 1<\alpha\le 2,\\[0.2cm]
\displaystyle y(a)=y(b)=y''(a)=0,
 \end{array}
 \right.
 \end{equation}
where  $D^{\alpha}$ is the Riemann-Liouville fractional derivative of order $\alpha,$    $q: [a,b]\to {\mathbb R}^+$ is a Lebesgue integrable  function and $f: {\mathbb R}^+\to {\mathbb R}^+$ is continuous.

The function $y \in C([a, b], {\mathbb R})$ is a solution to the problem (\ref{L-N-4}) if and only if $y$ satisfies the integral
equation
$$y(t) =\int_a^b G(t, s)q(s)f(y(s))ds,$$
where $G(t,s)$ is the Green's function defined by
\begin{equation}\label{Gr-L-N-3}
G(t,s)=\frac{1}{\Gamma(\alpha)}
\begin{cases}
\displaystyle\frac{(b-s)^{\alpha-1}(t-a)}{(b-a)}, &a\leq t\leq s\leq b,\\[0.3cm]
\displaystyle\frac{(b-s)^{\alpha-1}(t-a)}{(b-a)}-(t-s)^{\alpha-1}, &a\leq s \leq t \leq b.
\end{cases}
\end{equation}
The above  Green's function satisfies the following properties:
\begin{itemize}
\item[(i)] $G(t,s)\ge 0, $ for all $(t,s)\in [a,b]\times [a,b];$
\item[(ii)]~ For any $s\in [a,b],$
$$\max_{t\in [a,b]}G(t,s)=G(t_0,s)=\frac{(s-a)(b-s)^{\alpha-1}}{\Gamma(\alpha)(b-a)}+\frac{(\alpha-2)(b-s)^{\frac{(\alpha-1)^2}{\alpha-2}}}{(\alpha-1)^{\frac{\alpha-1}{\alpha-2}}(b-a)^{\frac{\alpha-1}{\alpha-2}}\Gamma(\alpha)},$$
where $\displaystyle t_0=s+\Big(\frac{(b-s)^{\alpha-1}}{(b-a)(\alpha-1)}\Big)^{\frac{1}{\alpha-2}}\in [s,b];$
\item[(iii)]~ $\displaystyle\max_{s\in [a,b]}G(t_0,s)\le\frac{(b-a)^{\alpha-1}}{\Gamma(\alpha)};$
\item[(iv)]~ $\displaystyle G(t,s)\ge \frac{(t-a)(b-t)}{(b-a)^2}G(t_0,s)$ for all $a\le t,s\le b.$
\end{itemize}

The following   Lyapunov-type inequalities are given in \cite{Ma}.
\begin{theorem}\label{t-N-4}
Assume that $f$ is bounded by two lines, that is, there exist two positive constants
$M$ and $N$ such that $Ny \le  f (y) \le My$ for any $y\in {\mathbb R}^+.$ If (\ref{L-N-4}) has a solution in $E^+=\{y\in C[a,b], y(t)\ge 0,~\mbox{for any}~ t\in [a,b]~ \mbox{and}~ \|y\|\ne 0\},$ then the
following   Lyapunov-type inequalities   hold:
\begin{itemize}
\item[(i)] $\displaystyle \int_a^b q(s)ds>\frac{\Gamma(\alpha)}{M(b-a)^{\alpha-1}};$
\item[(ii)]~ $\displaystyle \int_a^b(s-a)^2 (b-s)^{\alpha}q(s)ds\le\frac{4\Gamma(\alpha)(b-a)^3}{N};$
\item[(iii)]~ $\displaystyle \int_a^b(s-a)(b-s)^{\frac{\alpha^2-\alpha-1}{\alpha-2}}q(s)ds\le\frac{4\Gamma(\alpha)(b-a)^{\frac{3\alpha-2}{\alpha-2}}(\alpha-1)^{\frac{\alpha-1}{\alpha-2}}}{(\alpha-2)N}.$
\end{itemize}
\end{theorem}
The applications of the above inequalities are given in the following corollaries.
\begin{corollary}
For any $\displaystyle \lambda\in [0,\Gamma(\nu)]\cup\Big(\frac{4\Gamma(\nu)}{B(3,\nu+1)},+\infty\Big),$  where $\displaystyle B(x,y)=\int_0^1 s^{x-1}(1-s)^{y-1}ds, ~x>0, y>0,$ the
eigenvalue for the problem
 \begin{equation}\label{exam-9}
 \left\{\begin{array}{ll}
 (D^{\nu} y)(t)+\lambda y(t)=0,  \quad 0<t<1, ~ 2<\nu\le 3, \\
\displaystyle y(0)=y(1)=y''(0)=0,
 \end{array}
 \right.
 \end{equation}
has no corresponding eigenfunction $y\in E^+.$
\end{corollary}
In the next Corollary  we obtain an interval in
which the Mittag-Leffler function $E_{\nu,2}(z)$ with $\beta=2, 2<\nu\le 3$ has no real zeros.
\begin{corollary}
Let $2<\nu\le 3.$ Then the Mittag-Leffler function  {\small $\displaystyle E_{\nu,2}(z)=\sum_{k=0}^{\infty}\frac{z^k}{\Gamma(k\nu+2)}$}  has no real zeros for $\displaystyle z\in \Big(-\infty, -\frac{4\Gamma(\nu)}{B(3,\nu+1},\Big)\cup[-\Gamma(\nu), +\infty).$
\end{corollary}

In 2017,  Ru {\em et al.} \cite{RWAA}  obtained the Lyapunov-type inequality for the following {\em  fractional
Sturm-Liouville boundary value problem}
  \begin{equation}\label{Oct-3}
 \left\{\begin{array}{ll}
D^{\alpha}_{0^{+}}(p(t)y'(t))+q(t)y(t)=0, \quad 0<t<1, ~ 1<\alpha\le 2,\\
\displaystyle ay(0)-bp(0)y'(0)=0,~~ cy(1)+dp(1)y'(1)=0,
 \end{array}
 \right.
 \end{equation}
 where $a, b, c, d>0$ $D^{\alpha}$ is  the standard Riemann-Liouville fractional derivative of  order $\alpha,$ $p: [0,1]\to (0,\infty)$
  and  $q: [0,1]\to {\mathbb R}$ is a nontrivial Lebesgue  integrable function.

 The solution of the boundary value problem (\ref{Oct-3}) in terms of
 the integral equation is
\begin{equation*}
y(t)=\int_a^bG(t,s)q(s)y(s)ds,
\end{equation*}
where $G(t,s)$ is the Green's function given by
{\small $$G(t,s)=\frac{1}{\rho\Gamma(\alpha)}
\begin{cases}
\displaystyle \Big[b+a\int_0^t\frac{d\tau}{p(\tau)}\Big]\Big[d(1-s)^{\alpha-1}+c\int_s^1\frac{(\tau-s)^{\alpha-1}d\tau}{p(\tau)}\Big], &0\leq t\leq s\leq 1,\\[0.4cm]
\displaystyle  \Big[b+a\int_0^t\frac{d\tau}{p(\tau)}\Big]\Big[d(1-s)^{\alpha-1}+c\int_s^1\frac{(\tau-s)^{\alpha-1}d\tau}{p(\tau)}\Big]-H(t,s), &0\leq s \leq t \leq 1,
\end{cases}
$$}
$$\rho=bc+ac\int_0^1\frac{1}{p(\tau)}d\tau+ad, ~~ H(t,s)=a\Big[d+c\int_t^1\frac{d\tau}{p(\tau)}\Big]\int_0^t\frac{(\tau-s)^{\alpha-1}}{p(\tau)}d\tau.$$
Further, the above Green's function $G(t,s)$  satisfies the following properties:
\begin{itemize}
\item[(i)]~ $G(t,s)\ge 0$ for $0\le t,s\le 1;$
\item[(ii)] ~The   maximum value  of $G(t,s)$ is
\begin{equation}\label{over-G}
\overline G=\max_{0\le t,s\le 1}G(t,s)=\max\{\max_{s\in [0,1]}  G(s,s), \max_{s\in [0,1]}G(t_0(s),s)\},
\end{equation}
where
$$t_0(s)=s+\Bigg[\frac{ad(1-s)^{\alpha-1}+ac\int_s^1\frac{(\tau-s)^{\alpha-1}}{p(\tau)}d\tau}{\rho}\Bigg]^{\frac{1}{\alpha-1}}.$$
\end{itemize}

They obtained the following Lyapunov-type inequality for the problem  (\ref{Oct-3}).
\begin{theorem}\label{t-Oct-3}
For any nontrivial   solutions $y$ of the fractional boundary value problem (\ref{Oct-3}), the following  Lyapunov-type inequality holds:
$$\int_0^1|q(s)|ds>\frac{1}{\overline G},$$
where $\overline G$ is defined by (\ref{over-G}).
\end{theorem}

 In   \cite{RWAA},  the authors also considered the  {\em generalized  fractional
Sturm-Liouville boundary value problem}:
  \begin{equation}\label{Oct-3-f}
 \left\{\begin{array}{ll}
D^{\alpha}_{0^{+}}(p(t)y'(t))+q(t)f(y(t))=0, \quad 0<t<1, ~ 1<\alpha\le 2,\\[0.2cm]
\displaystyle ay(0)-bp(0)y'(0)=0,~~ cy(1)+dp(1)y'(1)=0,
 \end{array}
 \right.
 \end{equation}
 where $f: {\mathbb R}\to {\mathbb R}$ is a continuous function, and obtained the Lyapunov-type inequality for this problem as follows.

\begin{theorem}\label{t-Oct-3-f}
Let $f$ is a positive function on ${\mathbb R}.$ For any nontrivial   solutions $y$ of the fractional boundary value problem (\ref{Oct-3-f}), the following  Lyapunov-type inequality will be satisfied
$$\int_0^1|q(s)|ds>\frac{y^*}{\overline G \max_{y\in [y_*.y^*]}f(y)},$$
where $\overline G$ is defined by (\ref{over-G}) and $y_*=\min_{t\in [0,1]}y(t), y^*=\max_{t\in [0,1]}y(t).$
\end{theorem}

\section{Lyapunov inequalities for nonlocal   boundary value problems}

In 2017, Cabrera {\em et al.}  \cite{L7} obtained   Lyapunov-type inequalities for a {\em nonlocal} fractional
boundary value problem
  \begin{equation}\label{L-10}
 \left\{\begin{array}{ll}
^CD^{\alpha}y(t)+q(t)y(t)=0, \quad a<t<b,\\[0.3cm]
\displaystyle y'(a)=0,~~ \beta ^CD^{\alpha-1}y(b)+y(\eta)=0,
 \end{array}
 \right.
 \end{equation}
 where $^CD^{\alpha}$ is  the Caputo fractional derivative of fractional order $1<\alpha\le 2,$  $\beta>0, a < \eta <b,$ $\displaystyle \beta>\frac{(\beta-\eta)^{\alpha-1}}{\Gamma(\alpha)}$ and  $q: [a,b]\to {\mathbb R}$ is a continuous function.

 Note that  problem (\ref{L-10}) is the  fractional analogue of the classical nonlocal boundary value problem
 \begin{equation}\label{L-10-or}
 \left\{\begin{array}{ll}
 y''(t)+q(t)y(t)=0, \quad 0<t<1,\\[0.2cm]
\displaystyle y'(0)=0,~~ \beta y'(1)+y(\eta)=0, \, \, 0 < \eta <1,
 \end{array}
 \right.
 \end{equation}
which represents a thermostat model insulated at $t = 0$ with a controller dissipating heat at $t = 1$
depending on the temperature detected by a sensor at $t = \eta$  \cite{GuMe}.

 The function $y\in C([a,b], {\mathbb R})$ is a solution of the boundary value problem (\ref{L-10})
if and only if $y$ satisfies the integral equation
\begin{equation}\label{G-L-10}
y(t)=\int_a^bG(t,s)q(s)y(s)ds,
\end{equation}
where $G(t,s)$ is the Green's function defined by
\begin{equation}\label{Gr-L-10}
G(t,s)=\beta+H_{\eta}(s)-H_t(s),
\end{equation}
for $r\in [a,b], H_r:  [a,b]\to {\mathbb R}$ is
\begin{equation*}
H_r(s)=
\begin{cases}
\displaystyle \frac{(r-s)^{\alpha-1}}{\Gamma(\alpha)}, &a\leq s\leq r\leq b,\\
\displaystyle 0, &a\le r\le s  \leq b.
\end{cases}
\end{equation*}

The  Green's function defined in \eqref{Gr-L-10}   satisfies the relation
$$|G(t,s)|\le \beta+\frac{(\eta-a)^{\alpha-1}}{\Gamma(\alpha)}, ~~ (t,s)\in [a,b]\times [a,b].$$

 The Lyapunov-type inequality derived for the problem (\ref{L-10})  is  given in the following result.
 \begin{theorem}\label{L=10}
 If there exists a nontrivial continuous solution of the fractional
 boundary value problem (\ref{L-10}), then
 \begin{equation}\label{L10}
 \int_a^b  |q(s)|ds>\frac{\Gamma(\alpha)}{\beta \Gamma(\alpha)+(\eta-\alpha)^{\alpha-1}}.
 \end{equation}
 \end{theorem}
\vskip 0.3cm

In 2017, Cabrera {\em et al.} \cite{CSS} obtained a Lyapunov-type inequality for the following  {\em nonlocal fractional boundary value problem}
\begin{equation}\label{L-N-2}
 \left\{\begin{array}{ll}
(D^{\alpha}y)(t)+q(t)y(t)=0, \quad a<t<b,~~ 2<\alpha\le 3,\\[0.3cm]
\displaystyle y(a)=y'(a)=0, ~~ y'(b)=\beta y(\xi),
 \end{array}
 \right.
 \end{equation}
where  $D^{\alpha}$ is the Riemann-Liouville fractional derivative of order $\alpha,$ $a<\xi<b,$  $0\le \beta(\xi-a)^{\alpha-1}<(\alpha-1)(b-a)^{\alpha-2},$ and $q: [a,b]\to {\mathbb R}$ is a continuous function.

The unique solution of the nonlocal boundary value problem (\ref{L-N-2}) is given by
 {\small $$y(t)=\int_a^bG(t,s)q(s)y(s)ds+\frac{\beta(t-a)^{\alpha-1}}{(\alpha-1)(b-a)^{\alpha-2}-\beta(\xi-a)^{\alpha-1}}\int_a^bG(\xi,s)q(s)y(s)ds, $$}
where $G(t,s)$ is the Green's function defined by
\begin{equation}\label{Gr-L-N-2}
G(t,s)=\frac{1}{\Gamma(\alpha)}
\begin{cases}
\displaystyle\frac{(b-s)^{\alpha-2}(t-a)^{\alpha-1}}{(b-a)^{\alpha-2}}, &a\leq t\leq s\leq b,\\[0.3cm]
\displaystyle\frac{(b-s)^{\alpha-2}(t-a)^{\alpha-1}}{(b-a)^{\alpha-2}}-(t-s)^{\alpha-1}, &a\leq s \leq t \leq b.
\end{cases}
\end{equation}

The Green's function defined in (\ref{Gr-L-N-2}) satisfies the following properties:
\begin{itemize}
\item[(i)] $G(t,s)\ge 0, $ for all $(t,s)\in [a,b]\times [a,b];$
\item[(ii)]~ $G(t,s)$ is non-decreasing with respect to the first variable;
\item[(iii)]~ $0\le G(a,s)\le G(t,s)\le G(b,s),~ (t,s)\in [a,b]\times [a,b].$
\end{itemize}

Their Lyapunov-type inequality for the problem (\ref{L-N-2}) is expressed as follows.
\begin{theorem}\label{t-N-2}
If the problem (\ref{L-N-2}) has a nontrivial solution, then
 {\small \begin{equation}\label{LN2}
\int_a^b|q(s)|ds\ge\frac{\Gamma(\alpha)(\alpha-1)^{\alpha-1}}{(b-a)^{\alpha-1}(\alpha-2)^{\alpha-2}}\Big(1+\frac{\beta(b-a)^{\alpha-1}}{(\alpha-1)(b-a)^{\alpha-2}-\beta(\xi-a)^{\alpha-1}}\Big)^{-1}.
\end{equation}}
\end{theorem}

As an application of Theorem \ref{t-N-2}, we consider the following eigenvalue problem
\begin{equation}\label{exam-1}
 \left\{\begin{array}{ll}
D^{\alpha}y(t)+\lambda y(t)=0,  \quad a<t<b,~ 2<\alpha\le 3,\\[0.3cm]
\displaystyle y(a)=y'(a)=0, y'(b)=\beta y(\xi),
 \end{array}
 \right.
 \end{equation}
 where $a<\xi<b$ and $0\le \beta(\xi-a)^{\alpha-1}<(\alpha-1)(b-a)^{\alpha-2}.$ If $\lambda$ is an eigenvalue of problem (\ref{exam-1}), then
 $$|\lambda|\ge \frac{\alpha(\alpha-1)\Gamma(\alpha)}{(b-a)^{\alpha}}\Bigg(1+\frac{\beta(b-a)^{\alpha-1}}{(\alpha-1)(b-a)^{\alpha-2}-\beta(\xi-a)^{\alpha-1}}\Bigg)^{-1}.$$
This is an immediate consequence of Theorem \ref{t-N-2}.

\vskip 0.3cm

Very recently,  Y. Wang  and Q. Wang   \cite{WW}  obtained   Lyapunov-type inequalities for the fractional differential equations with {\em multi-point boundary conditions}
  \begin{equation}\label{Apr-1}
 \left\{\begin{array}{ll}
D^{\alpha} y(t)+q(t)y(t)=0, \quad a<t<b, ~ 2<\alpha\le3,\\
\displaystyle y(a)=y'(a)=0, ~~~ (D^{\beta+1}y)(b)=\sum_{i=1}^{m-2}b_i(D^{\beta}y)(\xi_i),
 \end{array}
 \right.
 \end{equation}
 where   $D^{\alpha}$ denotes the standard  Riemann-Liouville  fractional derivative of  order $\alpha,$  $\alpha>\beta+2,$ $0<\beta<1,$ $a<\xi_1<\xi_2<\ldots<\xi_{m-2}<b,$ $b_i\ge 0 (i=1,2,\ldots,m-2),$ $0\le \sum_{i=1}^{m-2}b_i(\xi_i-a)^{\alpha-\beta-1}<(\alpha-\beta-1)(b-a)^{\alpha-\beta-2}$
  and  $q: [a,b]\to {\mathbb R}$ is a continuous  function.

 The solution of the boundary value problem (\ref{Apr-1}) in terms of
 the integral equation is
\begin{equation*}
y(t)=\int_a^bG(t,s)q(s)y(s)ds+T(t)\int_a^b\Big(\sum_{i=1}^{m-2}b_iH(\xi,s)q(s)y(s)\Big)ds,,
\end{equation*}
where $G(t,s),$ $H(t,s)$ and $T(t)$ defined by
$$G(t,s)=\frac{1}{\Gamma(\alpha)}
\begin{cases}
\displaystyle \frac{(t-a)^{\alpha-1}(b-s)^{\alpha-\beta-2}}{(b-a)^{\alpha-\beta-2}}-(t-s)^{\alpha-1}, &a\leq s\leq t\leq b,\\[0.4cm]
\displaystyle  \frac{(t-a)^{\alpha-1}(b-s)^{\alpha-\beta-2}}{(b-a)^{\alpha-\beta-2}}, &a\leq t \leq s \leq b.
\end{cases}
$$
$$
H(t,s)=\frac{1}{\Gamma(\alpha)}
\begin{cases}
\displaystyle \frac{(t-a)^{\alpha-\beta-1}(b-s)^{\alpha-\beta-2}}{(b-a)^{\alpha-\beta-2}}-(t-s)^{\alpha-\beta-1}, &a\leq s\leq t\leq b,\\[0.4cm]
\displaystyle  \frac{(t-a)^{\alpha-\beta-1}(b-s)^{\alpha-\beta-2}}{(b-a)^{\alpha-\beta-2}}, &a\leq t \leq s \leq b.
\end{cases}
$$
$$T(t)=\frac{(t-a)^{\alpha-1}}{(\alpha-\beta-1)(b-a)^{\alpha-\beta-2}-\sum_{i=1}^{m-2}b_i(\xi_i-a)^{\alpha-\beta-1}}, ~~ t\ge a.$$
Further, the above   functions $G(t,s)$ and $H(t,s)$  satisfies the following properties:
\begin{itemize}
\item[(i)]~ $G(t,s)\ge 0$ for $a\le t,s\le b;$
\item[(ii)]~ $G(t,s)$ is non-decreasing with respect to the first variable;
\item[(iii)]~ $0\le G(a,s)\le G(t,s)\le G(b,s)=\frac{1}{\Gamma(\alpha)}(b-s)^{\alpha-\beta-2}[(b-a)^{\beta+1}-(b-s)^{\beta+1}],$ $ (t,s)\in [a,b]\times [a,b];$
\item[(iv)]~ for any $s\in [a,b],$
$$\max_{s\in [a,b]}G(b,s)=\frac{\beta+1}{\alpha-1}\Big(\frac{\alpha-\beta-2}{\alpha-1}\Big)^{(\alpha-\beta-2)/(\beta+1)}\frac{(b-a^{\alpha-1}}{\Gamma(\alpha)};$$
\item[(v)]~$H(t,s)\ge 0$ for $a\le t,s\le b;$
\item[(vi)]~~ $H(t,s)$ is non-decreasing with respect to the first variable;
\item[(vii)]~~~ $0\le H(a,s)\le H(t,s)\le H(b,s)=\frac{1}{\Gamma(\alpha)}(b-s)^{\alpha-\beta-2}(s-a), (t,s)\in [a,b]\times [a,b];$
\item[(viii)]~~~$$\max_{s\in [a,b]}H(b,s)=H(b,s^*)=\frac{(\alpha-\beta-2)^{\alpha-\beta-2}}{\Gamma(\alpha)}\Bigg(\frac{b-a}{\alpha-\beta-1}\Bigg)^{\alpha-\beta-1},$$
where $\displaystyle s^*=\frac{\alpha-\beta-2}{\alpha-\beta-1}a+\frac{1}{\alpha-\beta-1}b.$
\end{itemize}

They obtained the following Lyapunov-type inequalities.
\begin{theorem}\label{t-Apr-1}
If the fractional boundary value problem
 \begin{equation*}
 \left\{\begin{array}{ll}
D^{\alpha} y(t)+q(t)y(t)=0, \quad a<t<b, ~ 2<\alpha\le3,\\
\displaystyle y(a)=y'(a)=0, ~~~ (D^{\beta+1}y)(b)=\sum_{i=1}^{m-2}b_i(D^{\beta}y)(\xi_i),
 \end{array}
 \right.
 \end{equation*}
 has a nontrivial solution, where $q$ is a real and continuous function, then
 $$\int_a^b(b-s)^{\alpha-\beta-2}\Big[(b-a)^{\beta+1}-(b-s)^{\beta+1}+\sum_{i=1}^{m-2}b_iT(b)(s-a)\Big]|q(s)|ds\ge \Gamma(\alpha).$$
\end{theorem}

Note that
\begin{eqnarray*}
&&\Gamma(\alpha)\Big[G(b,s)+\sum_{i=1}^{m-2}b_iT(b)H(b,s)\Big]\\
&\le&\Gamma(\alpha)\Big[\max_{s\in [a,b]}G(b,s)+\sum_{i=1}^{m-2}b_iT(b)\max_{s\in [a,b]}H(b,s)\Big]\\
&=&\frac{\beta+1}{\alpha-1}\Bigg(\frac{\alpha-\beta-2}{\alpha-1}\Bigg)^{(\alpha-\beta-2)/(\beta+1)}(b-a)^{\alpha-1}\\
&&+\sum_{i=1}^{m-2}b_iT(b)(\alpha-\beta-2)^{\alpha-\beta-2}\Bigg(\frac{b-a}{\alpha-\beta-1}\Bigg)^{\alpha-\beta-1}.
\end{eqnarray*}
Thus we have
\begin{corollary}\label{Cor1-Apr-1}
If the fractional boundary value problem
 \begin{equation*}
 \left\{\begin{array}{ll}
D^{\alpha} y(t)+q(t)y(t)=0, \quad a<t<b, ~ 2<\alpha\le3,\\
\displaystyle y(a)=y'(a)=0, ~~~ (D^{\beta+1}y)(b)=\sum_{i=1}^{m-2}b_i(D^{\beta}y)(\xi_i),
 \end{array}
 \right.
 \end{equation*}
 has a nontrivial solution, where $q$ is a real and continuous function, then
 {\small \begin{eqnarray*}
 \int_a^b|q(s)|ds\ge\frac{\Gamma(\alpha)}{\frac{\beta+1}{\alpha-1}\Big(\frac{\alpha-\beta-2}{\alpha-1}\Big)^{\frac{\alpha-\beta-2}{\beta+1}}(b-a)^{\alpha-1}+\sum_{i=1}^{m-2}b_iT(b)(\alpha-\beta-2)^{\alpha-\beta-2}\Big(\frac{b-a}{\alpha-\beta-1}\Big)^{\alpha-\beta-1}}.
 \end{eqnarray*}}
\end{corollary}

If $\beta=0$  in Theorem \ref{t-Apr-1} we obtain
\begin{corollary}\label{Cor1-Apr-2}
If the fractional boundary value problem
 \begin{equation*}
 \left\{\begin{array}{ll}
D^{\alpha} y(t)+q(t)y(t)=0, \quad a<t<b, ~ 2<\alpha\le3,\\
\displaystyle y(a)=y'(a)=0, ~~~ y'(b)=\sum_{i=1}^{m-2}b_iy(\xi_i),
 \end{array}
 \right.
 \end{equation*}
 has a nontrivial solution, where $q$ is a real and continuous function, then
 \begin{eqnarray*}
 &&\int_a^b(b-s)^{\alpha-2}(s-a)|q(s)|ds\\
 &\ge&\frac{\Gamma(\alpha)}{1+\sum_{i=1}^{m-2}b_iT(b)}\\
 &=&\frac{(\alpha-\beta-1)(b-a)^{\alpha-\beta-2}+\sum_{i=1}^{m-2}b_i(\xi_i-a)^{\alpha-\beta-1}}{(\alpha-\beta-1)(b-a)^{\alpha-\beta-2}-\sum_{i=1}^{m-2}b_i(\xi_i-a)^{\alpha-\beta-1}+\sum_{i=1}^{m-2}b_i(b-a)^{\alpha-1}}\Gamma(\alpha).
 \end{eqnarray*}
\end{corollary}

If $\beta=0$ in Corollary \ref{Cor1-Apr-1} we have
\begin{corollary}\label{Cor-Apr-1}
If the fractional boundary value problem
 \begin{equation*}
 \left\{\begin{array}{ll}
D^{\alpha} y(t)+q(t)y(t)=0, \quad a<t<b, ~ 2<\alpha\le3,\\
\displaystyle y(a)=y'(a)=0, ~~~ y'(b)=\sum_{i=1}^{m-2}b_i y(\xi_i),
 \end{array}
 \right.
 \end{equation*}
 has a nontrivial solution, where $q$ is a real and continuous function, then
 \begin{eqnarray*}
&& \int_a^b|q(s)|ds\\
&\ge&\frac{\Gamma(\alpha)}{1+\sum_{i=1}^{m-2}b_iT(b)}\cdot \frac{(\alpha-1)^{\alpha-1}}{(b-a)^{\alpha-1}(\alpha-2)^{\alpha-2}}\\
&=&\frac{(\alpha-\beta-1)(b-a)^{\alpha-\beta-2}-\sum_{i=1}^{m-2}b_i(\xi_i-a)^{\alpha-\beta-1}}{(\alpha-\beta-1)(b-a)^{\alpha-\beta-2}-\sum_{i=1}^{m-2}b_i(\xi_i-a)^{\alpha-\beta-1}+\sum_{i=1}^{m-2}b_i(b-a)^{\alpha-1}}\\
&&\times \frac{\Gamma(\alpha)(\alpha-1)^{\alpha-1}}{(b-a)^{\alpha-1}(\alpha-2)^{\alpha-2}}.
 \end{eqnarray*}
\end{corollary}

\section{Lyapunov inequalities for fractional ${\bf p}$-Laplacian  boundary value problems}
 In this section we present  Lyapunov-type inequalities for  fractional $p$-Laplacian boundary value problems.

  In 2016, Al Arifi {\em et al.} \cite{AAJLS} considered the nonlinear fractional boundary value problem
\begin{equation}\label{p-L-1}
 \left\{\begin{array}{ll}
D^{\beta}(\Phi_p((D^{\alpha}y(t)))+q(t)\Phi_p(y(t))=0, \quad a<t<b,\\[0.3cm]
\displaystyle y(a)=y'(a)=y'(b)=0, ~~ D^{\alpha}y(a)=D^{\alpha}y(b)=0,
 \end{array}
 \right.
 \end{equation}
 where $2<\alpha\le 3,$ $1<\beta\le 2,$ $D^{\alpha}, D^{\beta}$ are the Riemann-Liouville fractional derivatives of orders $\alpha$ and $\beta$ respectively,  $\Phi_p(s)=|s|^{p-2}s, ~p>1$ is $p$-Laplacian operator  and $q:[a,b]\to {\mathbb R}$ is a continuous function.

For $h \in C([a, b], {\mathbb R})$, the  linear variant of the problem  (\ref{p-L-1}):
 \begin{equation}\label{p-L-1b}
 \left\{\begin{array}{ll}
D^{\beta}(\Phi_p((D^{\alpha}y(t)))+h(t)=0, \quad a<t<b,\\[0.3cm]
\displaystyle y(a)=y'(a)=y'(b)=0, ~~ D^{\alpha}y(a)=D^{\alpha}y(b)=0,
 \end{array}
 \right.
 \end{equation}
 has the unique solution
 $$y(t)=\int_a^b G(t,s)\Phi_q\Big(\int_a^b H(s,\tau)h(\tau) d\tau\Big)ds,$$
 where
 \begin{equation}\label{Gr-p-L-1b}
H(t,s)=\frac{1}{\Gamma(\beta)}
\begin{cases}
\displaystyle\Big(\frac{b-s}{b-a}\Big)^{\beta-1}(t-a)^{\beta-1}, &a\leq t\leq s\leq b,\\[0.4cm]
\displaystyle\Big(\frac{b-s}{b-a}\Big)^{\beta-1}(t-a)^{\beta-1}-(t-s)^{\beta-1}, &a\leq s \leq t \leq b,
\end{cases}
\end{equation}
  and $G(t,s)$ is the Green's function for the boundary value problem
\begin{equation}\label{p-L-1a}
 \left\{\begin{array}{ll}
D^{\beta}y(t)+h(t)=0, ~~2<\alpha\le 3,~~ a<t<b,\\
\displaystyle y(a)=y'(a)=y'(b)=0,  \end{array}
 \right.
 \end{equation}
which is   given by
\begin{equation}\label{Gr-p-L-2}
G(t,s)=\frac{1}{\Gamma(\alpha)}
\begin{cases}
\displaystyle\Big(\frac{b-s}{b-a}\Big)^{\alpha-2}(t-a)^{\alpha-1}, &a\leq t\leq s\leq b,\\[0.4cm]
\displaystyle\Big(\frac{b-s}{b-a}\Big)^{\alpha-2}(t-a)^{\alpha-1}-(t-s)^{\alpha-1}, &a\leq s \leq t \leq b.
\end{cases}
\end{equation}
Observe that the following estimates hold:
\begin{itemize}
\item[(i)] $0\le G(t,s)\le G(b,s), ~~ (t,s)\in [a,b]\times [a,b],$
\item[(ii)]~ $0\le H(t,s)\le H(s,s), ~~ (t,s)\in [a,b]\times [a,b].$
\end{itemize}
For the problem (\ref{p-L-1}),  the Lyapunov-type inequality is the following:
\begin{theorem}\label{t-p-L-1}
Let $2<\alpha\le 3,$ $1<\beta\le 2,$ $p>1,$    and $q\in C[a,b].$  If (\ref{p-L-1}) has a nontrivial solution, then
{\small\begin{equation}\label{p-eq-1}
\int_a^b(b-s)^{\beta-1}(s-a)^{\beta-1}|q(s)|ds\ge[\Gamma(\alpha)]^{p-1}\Gamma(\beta)(b-a)^{\beta-1}\Big(\int_a^b(b-s)^{\alpha-2}(s-a)ds\Big)^{1-p}.
\end{equation}}
\end{theorem}

 Now we  present an application  of this  result to eigenvalue problems.
 \begin{corollary}
 Let $\lambda$ be an eigenvalue of the problem
 \begin{equation}\label{exam-6}
 \left\{\begin{array}{ll}
D^{\beta}_{0}(\Phi_p(D^{\alpha}_{0+}y(t)))+\lambda \Phi_p(y(t))=0,  \quad 0<t<1, \\[0.3cm]
\displaystyle y(0)=y'(0)=y'(1)=0,~~ D^{\alpha}_{0+}y(0)=D^{\alpha}_{0+}y(1)=0,
 \end{array}
 \right.
 \end{equation}
 where $2<\alpha\le 3, 1<\beta\le 2, p>1,$ then
$$|\lambda|\ge \frac{\Gamma(2\beta)}{\Gamma(\beta)}\Bigg(\frac{\Gamma(\alpha)\Gamma(\alpha+1)}{\Gamma(\alpha-1)}\Bigg)^{p-1}.$$
 \end{corollary}
 In particular, for  $p = 2,$ that is, for  $\Phi_p(y(t))=y(t),$ the bound on $\lambda$ takes the form:
 $$|\lambda|\ge \frac{\Gamma(\alpha)\Gamma(\alpha+1)\Gamma(2\beta)}{\Gamma(\alpha-1)\Gamma(\beta)}.$$

In 2017, Liu {\em et al.} \cite{LXYB}  considered the nonlinear fractional $p$-Laplacian boundary value problem of the form:
\begin{equation}\label{p-L-2-x}
 \left\{\begin{array}{ll}
D^{\beta}(\Phi_p(^CD^{\alpha}y(t)))-q(t)f(y(t))=0, \quad a<t<b,\\
\displaystyle y'(a)={}^CD^{\alpha}y(a)=0, ~~ y(b)={}^CD^{\alpha}y(b)=0,
 \end{array}
 \right.
 \end{equation}
 where $1<\alpha, \beta\le 2,$   $^CD^{\alpha}, D^{\beta}$ are the Riemann-Liouville fractional derivatives of orders $\alpha$ and $\beta$ respectively,  $\Phi_p(s)=|s|^{p-2}s, p>1,$ and $q:[a,b]\to {\mathbb R}$ is a continuous function.

 An integral equation equivalent to the problem (\ref{p-L-2-x})
 is
 $$y(t)=\int_a^b G(t,s)\Phi_q\Big(\int_a^b H(s,\tau)f(y(\tau)) d\tau\Big)ds,$$
 where
 $$G(t,s)=\frac{1}{\Gamma(\alpha)}
\begin{cases}
\displaystyle(b-s)^{\alpha-1}, &a\leq t\leq s\leq b,\\[0.2cm]
\displaystyle (b-s)^{\alpha-1}-(t-s)^{\alpha-1}, &a\leq s \leq t \leq b,
\end{cases}
$$
and
 \begin{equation}\label{Gr-p-L-2-x}
H(t,s)=\frac{1}{\Gamma(\beta)}
\begin{cases}
\displaystyle\Big(\frac{s-a}{b-a}\Big)^{\beta-1}(b-s)^{\beta-1}, &a\leq t\leq s\leq b,\\[0.4cm]
\displaystyle\Big(\frac{s-a}{b-a}\Big)^{\beta-1}(b-s)^{\beta-1}-(t-s)^{\beta-1}, &a\leq s \leq t \leq b.
\end{cases}
\end{equation}
Moreover, the following estimates hold:
\begin{itemize}
\item[(i)] $  H(t,s)\ge 0$ for all $a\le t,s\le b;$
\item[(ii)]~ $\max_{t\in [a,b]}H(t,s)=H(s,s), ~ s\in [a,b];$
\item[(iii)]~ $H(t,s)$ has a unique maximum given by
$$\max_{s\in [a,b]}H(s,s)=\frac{(b-a)^{\beta-1}}{4^{\beta-1}\Gamma(\beta)};$$
\item[(iv)]~ $\displaystyle 0\le G(t,s)\le G(s,s)=\frac{1}{\Gamma(\alpha)}(b-s)^{\alpha-1}$ for all $a\le t,s\le b;$
\item[(v)]~ $G(t,s)$ has a unique maximum given by
$$\max_{s\in [a,b]}G(s,s)=\frac{1}{\Gamma(\alpha)}(b-a)^{\alpha-1}.$$
\end{itemize}

The Lyapunov-type inequalities for the problem (\ref{p-L-2-x}) are as follows.
\begin{theorem}\label{t-p-L-2-x}
Let $p:  [a,b]\to {\mathbb R}^+$ be a real   Lebesgue function. Suppose that there exists a positive constant  $M$ satisfying $0\le f(x)\le M \Phi_p(x)$ for any $x\in {\mathbb R}^+.$  If (\ref{p-L-2-x}) has a nontrivial solution in  $E^+=\{y\in C[a,b], y(t)\ge 0,~\mbox{for any}~ t\in [a,b]~ \mbox{and}~ \|y\|\ne 0\},$ then the following Lyapunov inequality holds:
\begin{equation}\label{t-L-2-x}
\int_a^b q(s)ds>\frac{4^{\beta-1}\Gamma(\beta)}{M(b-a)^{\beta-1}}\Phi_p\Big(\frac{\Gamma(\alpha+1)}{(b-a)^{\alpha}}\Big).
\end{equation}
\end{theorem}

\begin{theorem}\label{t-p-L-2-x-a}
Let $p:  [a,b]\to {\mathbb R}^+$ be a real   Lebesgue function. Assume that $f\in C({\mathbb R}^+, {\mathbb R}^+)$ is a concave and nondecreasing function. If (\ref{p-L-2-x}) has a nontrivial solution in  $E^+=\{y\in C[a,b], y(t)\ge 0,~\mbox{for any}~ t\in [a,b]~ \mbox{and}~ \|y\|\ne 0\},$ then the following Lyapunov inequality holds:
\begin{equation}\label{t-L-2-x-a}
\int_a^b q(s)ds>\frac{4^{\beta-1}\Gamma(\beta)\Phi_p(\Gamma(\alpha+1))\Phi_p(\eta)}{M(b-a)^{\alpha p+\beta-\alpha-1}f(\eta)},
\end{equation}
where $\eta=\max_{t\in [a,b]}y(t).$
\end{theorem}

 As an  application of the foregoing  results, we give the following corollary.
 \begin{corollary}
If $\lambda\in [0,4^{\beta-1}\Gamma(\beta)\Phi_p(\Gamma(\alpha+1))],$ then the following eigenvalue   problem
 \begin{equation}\label{exam-7}
 \left\{\begin{array}{ll}
D^{\beta}(\Phi_p(^CD^{\alpha}_{0+}y(t)))-\lambda \Phi_p(y(t))=0,  \quad 0<t<1, \\
\displaystyle y'(0)=^CD^{\alpha}=0,~~ y(1)=^CD^{\alpha}D^{\alpha}y(0)=0,
 \end{array}
 \right.
 \end{equation}
 has no corresponding eigenfunction $y\in E^+$, where $1<\alpha, \beta\le 2,$ and $p>1.$
 \end{corollary}

\section{Lyapunov inequalities for boundary value problems with  mixed fractional derivatives}
In 2017,  Guezane-Lakoud {\em et al.} \cite{LKT}  obtained a   Lyapunov-type  inequality for the following problem involving both right Caputo and left Riemann-Liouville fractional derivatives:
   \begin{equation}\label{e-x-1}
\left\{\begin{array} {ll}
\displaystyle -^CD^{\alpha}_{b-}D^{\beta}_{a+}y(t)+q(t)y(t),\  t\in [a,b],\\[0.3cm]
y(a)=D^{\beta}_{a+}y(b)=0,
\end{array} \right.
\end{equation}
where $0<\alpha, \beta<1, 1<\alpha+\beta\le 2,$ $^CD^{\alpha}_{b-}$ denotes the right Caputo fractional derivative, $D^{\beta}_{a+}$ denotes left Riemann-Liouville fractional derivative,  and   $q: [a,b]\to {\mathbb R}^+$ is a continuous  function.

The left and right Riemann-Liouville fractional integrals of order $p>0$ for a
function $g: (0,\infty)\rightarrow{\mathbb R}$ are  respectively defined  by
\begin{equation}\nonumber
I^{p}_{a+}g(t)=\int_a^t\frac{(t-s)^{p-1}}{\Gamma(p)}g(s)ds,
\end{equation}
\begin{equation}\nonumber
I^{p}_{b-}g(t)=\int_t^b\frac{(s-t)^{p-1}}{\Gamma(p)}g(s)ds,
\end{equation}
provided the right-hand sides are point-wise defined on $(0,\infty)$, where $\Gamma$ is the Gamma function.

The left Riemann-Liouville fractional derivative and the right Caputo fractional derivative of order $p>0$ for
a continuous function $g: (0,\infty)\to {\mathbb R}$ are respectively given by
\begin{equation}\nonumber
D^{p}_{a+}g(t)=\frac{d^n}{dt^n}(I^{n-p}_{a+})(t),
\end{equation}
\begin{equation}\nonumber
^cD^{p}_{b-}g(t)=(-1)^nI^{n-p}_{b-}g^{(n)}(t),
\end{equation}
 where $n-1<p<n.$

The function $y \in C[a, b]$ is a solution to the problem (\ref{e-x-1}) if and only if $y$ satisfies the integral
equation
$$y(t) =\int_a^b G(t, s)q(s)f(y(s))ds,$$
where $G(t,s)$ is the Green's function defined by
\begin{equation}\label{Gr-L-x-1}
G(t,s)=\frac{1}{\Gamma(\alpha)\Gamma(\beta)}
\begin{cases}
\displaystyle\int_a^r(t-s)^{\beta-1}(r-s)^{\alpha-12}, &a\leq r\leq t\leq b,\\[0.4cm]
\displaystyle\int_a^t(t-s)^{\beta-1}(r-s)^{\alpha-1}, &a\leq t \leq s \leq b.
\end{cases}
\end{equation}
The above  Green's function satisfies the following properties:
\begin{itemize}
\item[(i)] $G(t,s)\ge 0, $ for all $a\le r\le t\le b;$
\item[(ii)]~ $\max_{t\in [a,b]}G(t,r)=G(r,r)$ for all $r\in [a,b];$
\item[(iii)]~ $\displaystyle\max_{r\in [a,b]}G(r,r)=\frac{(b-a)^{\alpha+\beta-1}}{(\alpha+\beta-1)\Gamma(\alpha)\Gamma(\beta)}.$
\end{itemize}

The following result describes the Lyapunov inequality for problem (\ref{e-x-1}).
\begin{theorem}\label{t-x-1}
Assume that $0<\alpha, \beta<1$ and $1<\alpha+\beta\le 2.$ If the fractional
boundary value problem (\ref{e-x-1}) has a nontrivial continuous solution, then
$$\int_a^b|q(r)|dr\ge \frac{(\alpha+\beta-1)\Gamma(\alpha)\Gamma(\beta)}{(b-a)^{\alpha+\beta-1}}.$$
\end{theorem}

\section{Lyapunov inequality for Hadamard type fractional boundary value problems}
 Let us begin this section with some fundamental definitions.
\begin{definition} \cite{Kil} The Hadamard derivative of fractional order $q$ for a function $g: [1, \infty)$ $\to \mathbb{R}$ is
defined as
 $$^{H}D^q g(t)=\frac{1}{\Gamma(n-q)}\left(t\frac{d}{dt}\right)^n\int_{1}^t \left(\log\frac{t}{s}\right)^{n-q-1}\frac{g(s)}{s}ds, ~~n-1 < q < n,~
n=[q]+1,$$   where $[q]$ denotes the integer part of the real
 number $q$ and $\log (\cdot) =\log_e (\cdot)$ is the usual Napier logarithm.
 \end{definition}
\begin{definition}\cite{Kil} The Hadamard fractional integral of order $q$ for a function $g$ is defined as
$$I^q g(t)=\frac{1}{\Gamma(q)}\int_1^t \left(\log\frac{t}{s}\right)^{q-1}\frac{g(s)}{s}ds, ~~q>0, ~$$
provided the integral exists.
\end{definition}

 In 2017, Ma {\em et al.} \cite{L6}   obtained  a Lyapunov-type inequality for a  {\em Hadamard fractional
boundary value problem}
  \begin{equation}\label{L-6}
 \left\{\begin{array}{ll}
^HD^{\alpha}y(t)-q(t)y(t)=0, \quad 1<t<e,\\[0.2cm]
\displaystyle y(1)=y(e)=0,
 \end{array}
 \right.
 \end{equation}
 where $^HD^{\alpha}$ is the fractional derivative in the sense of the Hadamard of order $1<\alpha\le 2$  and  $q: [a,b]\to {\mathbb R}$ is a continuous function.

 The function  $y\in C([1,e], {\mathbb R})$ is a solution of the boundary value problem (\ref{L-6})
if  and only if  $y$ satisfies the integral equation
\begin{equation}\label{G-L-6}
y(t)=\int_a^bG(t,s)q(s)y(s)ds,
\end{equation}
where $G(t,s)$ is the Green's function defined by
 {\small \begin{equation}\label{Gr-L-6}
G(t,s)=\frac{1}{\Gamma(\alpha)}
\begin{cases}
\displaystyle -\Big(\log\frac{e}{s}\Big)^{\alpha-1}\frac{(\log t)^{\alpha-1}}{s}, &1\leq t\leq s\leq e,\\[0.4cm]
\displaystyle-\Big(\log\frac{e}{s}\Big)^{\alpha-1}\frac{(\log t)^{\alpha-1}}{s}+\Big(\log\frac{t}{s}\Big)^{\alpha-1}\frac{1}{s}, &1\leq s \leq t \leq e.
\end{cases}
\end{equation}}
such that
$$|G(t,s)|\le  \frac{\lambda^{\alpha-1}}{\Gamma(\alpha)} (1-\lambda)^{\alpha-1}\exp (-\lambda),$$
with
\begin{equation}\label{lam}
 \lambda=\frac{1}{2}\Big(2\alpha-1-\sqrt{(2\alpha-2)^2+1}~\Big).
 \end{equation}

 The result concerning the Lyapunov-type inequality for the problem (\ref{L-6}) is  as follows.
 \begin{theorem}\label{L=6}
 If there exists a nontrivial continuous solution of the fractional
 boundary value problem (\ref{L-6}), then
 \begin{equation}\label{L6}
 \int_a^b |q(s)|ds>\Gamma(\alpha) \lambda^{1-\alpha}.(1-\lambda)^{1-\alpha}\exp \lambda,
 \end{equation}
 where $ \lambda$ is defined by (\ref{lam}).
 \end{theorem}

 For recent results on Hadamard type fractional boundary value problems, we refer the interested reader to the book \cite{AANT-B2}.

\section{Lyapunov inequality for boundary value problems with  Prabhakar  fractional derivative}
In  \cite{AbPa},  the authors discussed  Lyapunov-type inequality for the following fractional boundary value problem
involving the  $k$-Prabhakar derivative:
\begin{equation}\label{p-L-3-x}
 \left\{\begin{array}{ll}
(_k D^{\gamma}_{\rho,\beta,\omega,a+}y)(t)+q(t)f(y(t))=0, \quad a<t<b,\\[0.2cm]
\displaystyle y(a)=y(b)=0,
 \end{array}
 \right.
 \end{equation}
where $_k D^{\gamma}_{\rho,\beta,\omega,a+}$ is the $k$-Prabhakar differential operator of order $\beta\in (1,2],$ $ k\in {\mathbb R}^+$ and $\rho, \gamma, \omega \in {\mathbb C}.$
 The $k$-Prabhakar integral operator is defined as
 $$(_k P_{\alpha, \beta, \omega}\phi)(t)=\int_0^x\frac{(x-t)^{\frac{\beta}{k}-1}}{k}E_{k,\alpha,\beta}[\omega(x-t)^{\frac{\alpha}{k}}]\phi(t)dt, ~~ x>0,$$
 where $E_{k,\alpha,\beta}$ is the $k$-Mittag-Leffler function given by
 $$E_{k,\alpha,\beta}^{\gamma}(z)=\sum_{n=0}^{\infty}\frac{(\gamma)_{n,k}z^n}{\Gamma_k(\alpha n+\beta)n!},$$
  $\Gamma_k(x)$ is the $k$-Gamma function $\displaystyle \Gamma_k(x)=\lim_{n\to \infty}\frac{n!k^n(nk)^{\frac{x}{k}-1}}{(x)_{n,k}}$ and $\displaystyle (\gamma)_{n,k}=\frac{\Gamma_k(\gamma+nk)}{\Gamma_k(\gamma)}$ is the Pochhammer $k$-symbol.

The $k$-Prabhakar derivative is defined as
$$_kD^{\gamma}_{\rho,\beta,\omega}f(x)=\Big(\frac{d}{dx}\Big)^mk^m_kP^{-\gamma}_{\rho,mk-\beta,\omega}f(x),$$
 where $m=[\beta/k]+1.$

  An integral equation related to the problem (\ref{p-L-3-x})
is
 $$y(t)=\int_a^b G(t,s)q(s)y(s)ds,$$
 where
 {\small  \begin{equation}\label{Prab}
 G(t,s)=
\begin{cases}
\displaystyle \frac{(t-a)^{\frac{\beta}{k}-1} E_{k,\rho,\beta}^{\gamma}(\omega(t-a)^{\frac{\rho}{k}})(b-s)^{\frac{\beta}{k}-1}}{(b-a)^{\frac{\beta}{k}-1} E_{k,\rho,\beta}^{\gamma}(\omega(b-a)^{\frac{\rho}{k}})k} E_{k,\rho,\beta}^{\gamma}(\omega(b-s)^{\frac{\rho}{k}}),&a\leq t\leq s\leq b,\\[0.4cm]
\displaystyle \frac{(t-a)^{\frac{\beta}{k}-1} E_{k,\rho,\beta}^{\gamma}(\omega(t-a)^{\frac{\rho}{k}})(b-s)^{\frac{\beta}{k}-1}}{(b-a)^{\frac{\beta}{k}-1} E_{k,\rho,\beta}^{\gamma}(\omega(b-a)^{\frac{\rho}{k}})k} E_{k,\rho,\beta}^{\gamma}(\omega(b-s)^{\frac{\rho}{k}}\\
-\frac{(t-s)^{\frac{\beta}{k}-1}}{k}E_{k,\rho,\beta}^{\gamma}(\omega(t-s)^{\frac{\rho}{k}}), &a\leq s \leq t \leq b,
\end{cases}
\end{equation}}
which satisfies the following properties:
\begin{itemize}
\item[(i)] $G(t,s)\ge 0$ for all $a\le t,s\le b;$
\item[(ii)]~ $\max_{t\in [a,b]}G(t,s)=G(s,s)$ for all $s\in [a,b];$
\item[(iii)]~ $G(t,s)$ has a unique maximum given by
$$\max_{s\in [a,b]}G(s,s)=\Big(\frac{b-a}{4}\Big)^{\frac{\beta}{k}-1}\frac{E_{k,\rho,\beta}^{\gamma}\Big(\omega\Big(\frac{b-a}{2}\Big)^{\frac{\rho}{k}}\Big)E_{k,\rho,\beta}^{\gamma}\Big(\omega\Big(\frac{b-a}{2}\Big)^{\frac{\rho}{k}}\Big)}{kE_{k,\rho,\beta}^{\gamma}\Big(\omega (b-a)^{\frac{\rho}{k}}\Big)}.$$
\end{itemize}

The Lyapunov-type inequality for the problem (\ref{p-L-3-x})  is given in the  following result.
\begin{theorem}
If the problem (\ref{p-L-3-x}) has a nontrivial solution, then
$$\int_a^b|q(s)|ds\ge \Big(\frac{4}{b-a}\Big)^{\frac{\beta}{k}-1}\frac{kE_{k,\rho,\beta}^{\gamma}\Big(\omega (b-a)^{\frac{\rho}{k}}\Big)}{E_{k,\rho,\beta}^{\gamma}\Big(\omega\Big(\frac{b-a}{2}\Big)^{\frac{\rho}{k}}\Big)E_{k,\rho,\beta}^{\gamma}\Big(\omega\Big(\frac{b-a}{2}\Big)^{\frac{\rho}{k}}\Big)}.$$
\end{theorem}

The special case $k=1$ for the  problem (\ref{p-L-3-x}) has recently been studied in \cite{EsAn}.
\vskip 0.3cm
In 2017, Pachpatte  {\em et al.} \cite{PAK}  established some Lyapunov-type inequalities for the following {\em hybrid fractional    boundary value problem}
\begin{equation}\label{exam-16}
\begin{cases}
\displaystyle D^{\gamma}_{\rho, \beta, \omega}\left[\displaystyle\frac{y(t)}{f(t,y(t))}-\displaystyle\sum_{i=1}^n E^{\gamma}_{\rho, \beta, \omega}h_i(t,y(t))\right]+g(t)y(t)=0, \quad t\in(a,b),\\[4mm]
y(a)=y(b)=0,
\end{cases}
\end{equation}
where $D^{\gamma}_{\rho, \beta, \omega}$ denotes the Prabhakar fractional derivative of order  $\beta\in(1,2]$ starting from a point $a$,  $y\in C([a,b], \mathbb{R})$, $g\in L^1((a,b], \mathbb{R})$, $f\in C^1([a,b]\times\mathbb{R}, \mathbb{R}\setminus\{0\})$, $h_i\in C([a,b]\times\mathbb{R}, \mathbb{R})$, $\forall i=1,2,\ldots,n$ and  $E^{\alpha}_{\rho, \mu, \omega}$ is the Prabhakar fractional integral of order $\mu$ with the lower limit at the point $a$.
\vskip 0.3cm
We consider two cases:  (I) $h_i=0, i=1,2,\ldots,n$ and (II) $h_i\ne 0, i=1,2,\ldots,n.$\\
{\bf Case I: $h_i=0, i=1,2,\ldots,n.$ }
We consider the problem (\ref{exam-16})  with $h_i(t,\cdot)=0$ for all $t\in [a,b]$. For $\alpha\in(2,3]$, we first construct a Green's function for the following boundary value problem
\begin{equation}\label{exam-16a}
\begin{cases}
D^{\alpha}_{\rho, \mu, \omega}\left[\displaystyle\frac{y(t)}{f(t,y(t))}\right]+g(t)y(t)=0,\quad t\in(a,b), \\[4mm]
y(a)=y(b)=0,
\end{cases}
\end{equation}
with the assumption that $f$ is continuously differentiable and $f(t,y(t))\neq 0$ for all $t\in[a,b]$.
Let $y\in AC([a,b],\mathbb{R})$ be a solution of the problem \eqref{exam-16a}. Then the function $y$   satisfies the following integral equation
\begin{equation}\label{greensol-16a}
y=f(t,y)\int_a^bG(t,s)g(s)y(s)ds,
\end{equation}
where $G(t,s)$ is the Green's function defined by (\ref{Prab}).
 The Lyapunov-type inequality for this case is as follows.
\begin{theorem}
If the problem (\ref{exam-16a}) has a nontrivial solution, then
$$\int_a^b|q(s)|ds\ge \frac{1}{\|f\|}\Big(\frac{4}{b-a}\Big)^{\beta-1}\frac{E_{\rho,\beta}^{\gamma}\Big(\omega (b-a)^{\rho}\Big)}{E_{\rho,\beta}^{\gamma}\Big(\omega\Big(\frac{b-a}{2}\Big)^{\rho}\Big)E_{\rho,\beta}^{\gamma}\Big(\omega\Big(\frac{b-a}{2}\Big)^{\rho}\Big)},$$
where $\|f\|=\displaystyle\sup_{t\in[a,b], y\in {\mathbb R}}   |f(t,y)| $.
\end{theorem}
  {\bf Case II: $h_i\ne 0, i=1,2,\ldots, n.$}
Let $y\in AC[a,b]$ be a solution of the problem (\ref{exam-16}) given by
\begin{equation}
y(t)=f(t,y(t)) \int_a^b G(t,s)g\Big[(s)y(s)-\sum_{i=1}^n\int_a^bh_i(s,y(s))\Big]ds,
\end{equation}
where   $G(t,s)$  is defined as in \eqref{Prab}.

\begin{theorem}\label{thmp201-Pr} (Lyapunov-type inequality)
Assume that $|q(t)y(t)-\sum_{i=1}^{n}h_i(t,y(t))|\leq K|qt)|\|y\|, ~ K\in {\mathbb R}.$  If a nontrivial solution for the problem \eqref{exam-16} exists, then
$$\int_a^b|q(s)|ds\ge \frac{1}{K\|f\|}\Big(\frac{4}{b-a}\Big)^{\beta-1}\frac{E_{\rho,\beta}^{\gamma}\Big(\omega (b-a)^{\rho}\Big)}{E_{\rho,\beta}^{\gamma}\Big(\omega\Big(\frac{b-a}{2}\Big)^{\rho}\Big)E_{\rho,\beta}^{\gamma}\Big(\omega\Big(\frac{b-a}{2}\Big)^{\rho}\Big)}.$$
\end{theorem}

\section{Lyapunov inequality for $q$-difference  boundary value problems}

Let a $q$-real number denoted by $[a]_q$ be defined by
$$[a]_q=\frac{1-q^a}{1-q}, \, \, \, a \in \mathbb{R}, \, \, q \in \mathbb{R}^+\setminus \{1\}. $$

The $q$-analogue of the Pochhammer symbol ($q$-shifted factorial)
is defined as
$$(a; q)_{0}=1, ~~ (a; q)_{k}=\prod_{i=0}^{k-1}(1-aq^{i}), ~ k\in {\mathbb N} \cup \{\infty\}.$$
The $q$-analogue of the exponent $(x-y)^k$ is
$$(x-y)^{(0)}_a=1, ~~ (x-y)^{(k)}_a=\prod_{j=0}^{k-1}((x-a)-(y-a)q^{j}), ~ k\in {\mathbb N}, ~ x, y \in \mathbb{R}.$$

More generally, if $\gamma\in {\mathbb R},$ then
$$(x-y)^{(\gamma)}_a=(x-a)^{\gamma}\prod_{i=0}^{\infty}\frac{(x-a)-q^i(y-a)}{(x-a)-q^{\gamma+1}(y-a)}.$$
The $q$-Gamma function $\Gamma_q(y)$ is defined as
$$\Gamma_q(y)=\frac{(1-q)^{(y-1)}_0}{(1-q)^{y-1}}, $$
where $y \in \mathbb{R} \setminus \{0, -1, -2, \ldots\}.$ Observe that $\Gamma_q(y+1)=[y]_q \Gamma_q(y).$

The $q$-derivative of a function $f : [a; b]\to {\mathbb R} (a < b)$ is defined by
$$(_aD_qf)(t)=\frac{f(t)-f(qt+(1-q)a}{(1-q)(t-a)},~~ t\ne a$$
and
$$(_aD_qf)(a)=\lim_{t\to a}(_aD_qf)(t).$$

In 2016, Jleli and Samet \cite{MJBS} established a Lyapunov-type inequality for a Dirichlet boundary value problem of  fractional $q$-difference equations given by
\begin{equation}\label{p-L-5-x}
 \left\{\begin{array}{ll}
(_aD^{\alpha}y)(t)+\phi(t)y(t)=0,  \quad a<t<b,~ q\in [0,1),~ 1<\alpha\le 2,\\[0.2cm]
\displaystyle y(a)=y(b)=0,
 \end{array}
 \right.
 \end{equation}
where $_aD^{\alpha}$ denotes the fractional $q$-derivative of Riemann-Liouville type and $\phi: [a, b]\to {\mathbb R}$ is a continuous
function.

The solution $y\in C([a,b], {\mathbb R})$ of the problem (\ref{p-L-5-x}) satisfies the integral equation
$$y(t)=\int_a^tG_1(t, qs+(1-q)a)\phi(s)y(s)~_ad_qs+\int_t^bG_2(t,s)\phi(s)y(s)~_ad_qs, ~~ a\le t\le b,$$
where
$$G_1(t,s)=\frac{1}{\Gamma_q(\alpha)}\Big(\frac{(t-a)^{\alpha-1}}{(b-a)^{\alpha-1}}(b-s)_a^{(\alpha-1)}-(t-s)_a^{(\alpha-1)}\Big), ~a\le s\le t\le b,$$
 $$G_2(t,s)=\frac{(t-a)^{\alpha-1}}{\Gamma_q(\beta)(b-a)^{\alpha-1}}(b-(qs+(1-q)a)_a^{(\alpha-1)},~ a\le t\le s\le b,$$
 satisfying the relations
 \begin{itemize}
 \item[(i)] $0\le G_1(t,qs+(1-q)a)\le G_2(s,s), ~~ a<s\le t\le b;$
  \item[(ii)]~  $G_2(a,s)=0\le G_2(t,s)\le G_2(s,s), ~~ a\le t\le s\le b.$
  \end{itemize}

\begin{theorem} (Lyapunov-type inequality)
If the problem (\ref{p-L-5-x}) has a nontrivial solution, then
$$\int_a^b (s-a)^{\alpha-1}(b-(qs+(1-q)a))_a^{(\alpha-1)}|\phi(s)|~ _ad_qs\ge \Gamma_q(\alpha)(b-a)^{\alpha-1}.$$
\end{theorem}

 Taking $\alpha=2$ in the above theorem we have the following corollary.
 \begin{corollary}
 If a nontrivial continuous solution to the $q$-difference boundary value problem,
 \begin{equation}\label{p-L-5-x-c}
 \left\{\begin{array}{ll}
(_aD^{2}_qy)(t)+\phi(t)y(t)=0,  \quad a<t<b,~ q\in (0,1),\\[0.2cm]
\displaystyle y(a)=y(b)=0,
 \end{array}
 \right.
 \end{equation}
 exists, where $\phi: [a,b]\to {\mathbb R}$ is a continuous function, then
 $$\int_a^b(s-a)(b-(qs+(1-q)a))|\phi(s)|~ _ad_qs\ge (b-a).$$
 \end{corollary}

Some recent work on $q$-difference  boundary value problems can be found in  \cite{ANT-B1}.

\section{Lyapunov inequality for   boundary value problems involving a fractional derivative with respect to a certain function}

In 2017, Jleli {\em et al.} \cite{LKS-0}  considered the following fractional boundary value problem involving a fractional derivative with respect to a certain function $g$
 \begin{equation}\label{p-L-6-x}
 \left\{\begin{array}{ll}
(D^{\alpha}_{a,g}y)(t)+q(t)y(t)=0,  \quad a<t<b,~ 1<\alpha\le 2,\\
\displaystyle y(a)=y(b)=0,
 \end{array}
 \right.
 \end{equation}
 where $D^{\alpha}_{a,g}$ is the fractional derivative operator of order $\alpha$ with respect to a  nondecreasing function $g\in C^1([a,b], {\mathbb R})$ with $g'(x)>0,$
 for all $x\in [a,b],$ and $q:[a,b]\to {\mathbb R}$ is a continuous function.

 Let $f\in L^1((a,b), {\mathbb R}).$ The fractional integral of order $\alpha>0$ of $f$ with respect to the function $g$ is defined by
 $$(I_{a,g}^{\alpha} f)(t)=\frac{1}{\Gamma(\alpha)}\int_a^t\frac{g'(s)f(s)}{(g(t)-g(s))^{1-\alpha}}ds, ~ \mbox{a.e.}~ t\in [a,b].$$
Let $\alpha>0$ and $n$ be the smallest integer
greater than or equal to $\alpha.$ Let  $f:[a,b]\to {\mathbb R}$ be a function
such that   $((1/g'(t))(d/dt)^nI^{n-\alpha}_{a,g}f$  exists almost everywhere
on $[a,b].$  Then the fractional derivative of order $\alpha$ of $f$
with respect to the function $g$ is defined as
 \begin{eqnarray*}
 D^{\alpha}_{a,g}f(t)&=&\Big(\frac{1}{g'(t)}\frac{t}{dt}\Big)^nI^{n-\alpha}_{a,g}f(t)\\
 &=&\frac{1}{\Gamma(n-\alpha)}\Big(\frac{1}{g'(t)}\frac{t}{dt}\Big)^n\int_a^t\frac{g'(s)f(s)}{(g(t)-g(s))^{\alpha-n+1}}ds, ~\mbox{for}~ \mbox{a.e}~ t\in [a,b].
 \end{eqnarray*}

 \begin{theorem}\label{t-g}
 Assume that $q\in C([a,b], {\mathbb R})$  and $g\in C^1([a,b], {\mathbb R})$ be a nondecreasing function with $g'(x)>0,$ for all $x\in [a,b].$ If the problem  (\ref{p-L-6-x}) has a nontrivial solution,  then
 \begin{equation}\label{eq-g}
 \int_a^b[(g(s)-g(a)((g(b)-g(s))]^{\alpha-1}g'(s)|q(s)|ds\ge \Gamma(\alpha)(g(b)-g(a))^{\alpha-1}.
 \end{equation}
 \end{theorem}
 From inequality (\ref{eq-g}), we can  obtain
Lyapunov-type inequalities for different choices of the function $g.$ For instance, for $g(x)=x^{\beta}, x\in [a,b]$ and $g(x)=\log x, x\in [a,b]$ we have respectively  the following results.
 \begin{corollary}
  If the problem (\ref{p-L-6-x}) has a nontrivial solution and
 $g(x)=x^{\beta}, x\in [a,b], 0<a<b,$ then
 $$\int_a^b|q(s)|ds\ge \frac{\Gamma(\alpha)(b^{\beta}-a^{\beta})^{\alpha-1}}{\beta\phi_{\alpha,\beta}(s^*(\alpha,\beta))},$$
 where $\phi_{\alpha,\beta}(s^*(\alpha,\beta))=\max\{\phi_{\alpha,\beta}(s): s\in [a,b]\}>0.$
 \end{corollary}
 Taking $g(x)=\log x, x\in [a,b], 0<a<b,$ in Theorem \ref{t-g}, we deduce the
following Hartman-Wintner-type inequality
 $$\int_a^b\Big[\Big(\log\frac{s}{a}\Big)\Big(\log\frac{b}{s}\Big)\Big]^{\alpha-1}\frac{|q(s)|}{s}ds\ge \Gamma(\alpha)\Big(\log\frac{b}{a}\Big)^{\alpha-1},$$
for the Hadamard fractional boundary value problem of the form:
  \begin{equation}\label{p-L-6-x-H}
 \left\{\begin{array}{ll}
(^HD^{\alpha}_{a}y)(t)+q(t)y(t)=0,  \quad a<t<b,~ 1<\alpha\le 2,\\[0.2cm]
\displaystyle y(a)=y(b)=0.
 \end{array}
 \right.
 \end{equation}

\vskip 0.3cm
 In order to demonstrate the application of Theorem \ref{t-g}, we consider the eigenvalue problem
\begin{equation}\label{exam-5}
 \left\{\begin{array}{ll}
D^{\alpha}_{a,g}y(t)+\lambda y(t)=0,  \quad a<t<b, \\
\displaystyle y(a)=y(b)=0,
 \end{array}
 \right.
 \end{equation}
  and use the Lyapunov-type inequality (\ref{eq-g}) to obtain the following result.
 \begin{theorem} If $\lambda$ is an eigenvalue of fractional boundary value problem (\ref{exam-5}), then the following inequality holds
 $$|\lambda|\ge \frac{\Gamma(\alpha)(g(b)-g(a))^{\alpha-1}}{\int_{g(a)}^{g(b)} (x-g(a))^{\alpha-1}(g(b)-x)^{\alpha-1} dx}.$$
 \end{theorem}

\section{Lyapunov inequality for   boundary value problems involving left and right derivatives}
The left and right Caputo fractional derivatives are defined via  the Riemann-Liouville fractional derivatives (see \cite[p. 91]{Kil}). In particular, they are defined for a class of absolutely continuous functions.
\vskip 0.5cm
\begin{definition}  {\rm (Left and Right Riemann-Liouville Fractional Integrals \cite{Kil})}. Let $f$ be a function defined on $[a, b].$ The left and right Riemann-Liouville fractional integrals of order $\gamma$ for function $f$ denoted by $I^{\gamma}_{a+}$ and $I^{\gamma}_{b-}$  respectively, are defined by
$$I^{\gamma}_{a+}f(t)=\frac{1}{\Gamma(\gamma)}\int_a^t(t-s)^{\gamma-1}f(s)ds,~ t\in [a,b], ~ \gamma>0,$$
and
$$I^{\gamma}_{b-}f(t)=\frac{1}{\Gamma(\gamma)}\int_t^b(t-s)^{\gamma-1}f(s)ds,~ t\in [a,b], ~ \gamma>0,$$
provided the right-hand sides are pointwise defined on $[a, b],$ where $\Gamma > 0$ is the gamma function.
\end{definition}
\begin{definition}  {\rm (Left and Right Riemann-Liouville Fractional Derivatives \cite{Kil})}. Let $f$ be a function defined on $[a, b].$ The left
and right Riemann-Liouville fractional derivatives of order $\gamma$ for function $f$ denoted by $_aD^{\gamma}_t f (t)$ and $_tD^{\gamma}_b f (t),$  respectively,
are defined by
$$_aD^{\gamma}_t f (t)=\frac{d^n}{dt^n}_aD^{\gamma-n}_t f (t)=\frac{1}{\Gamma(n-\gamma)}\frac{d^n}{dt^n}\Big(\int_a^t(t-s)^{n-\gamma-1}f(s)ds\Big)$$
and
$$_tD^{\gamma}_b f (t)=(-1)^n\frac{d^n}{dt^n}_tD^{\gamma-n}_b f (t)=(-1)^n\frac{1}{\Gamma(n-\gamma)}\frac{d^n}{dt^n}\Big(\int_t^b(t-s)^{n-\gamma-1}f(s)ds\Big)$$
where $t \in  [a, b],$ $ n - 1 \le \gamma < n$ and $n\in{\mathbb N}.$
\end{definition}
\vskip 0.3cm

In 2017, Chen {\em et al.} \cite{CLL} obtained a Lyapunov-type inequality for the following boundary value problem
\begin{equation}\label{M-21}
 \left\{\begin{array}{ll}
\displaystyle  \frac{d}{dt}\Big(\frac{1}{2} I^{\beta}_{a+}y'(t)+\frac{1}{2}I^{\beta}_{b-}y'(t)\Big)+q(t)y(t)=0,~~ a<t<b,\\[0.3cm]
\displaystyle y(a)=y(b)=0,
 \end{array}
 \right.
 \end{equation}
where $I^{\beta}_{a+}$ and $I^{\beta}_{b-}$  denote the left and right Riemann-Liouville fractional integrals of order
$0\le \beta<1,$ respectively, and $q\in L^1([a,b], {\mathbb R}).$

\begin{theorem}
Let $q\in L^1([a,b], {\mathbb R})$ be a nonnegative function and there exists  a nontrivial solution for the boundary value problem (\ref{M-21}). Then
$$\int_a^b|q(s)|dt\ge\Bigg(\frac{2(b-a)^{\alpha-1/2}}{\Gamma(\alpha)(2\alpha-1)^{1/2}}\frac{1}{|\cos(\pi\alpha)|^{1/2}}\Bigg)^{-2}, ~~ \alpha=1-\frac{\beta}{2}.$$
\end{theorem}
\vskip 0.3cm
Jleli {\em et al.} \cite{JlKiSa-LR} considered the following quasilinear problem involving both left and right Riemann-Liouville fractional derivative operators:
\begin{equation}\label{LR-1}
 \left\{\begin{array}{ll}
\displaystyle \frac{1}{2}\Big(\, _tD^{\alpha}_b(|\,_aD^{\alpha}_t|^{p-2}_aD^{\alpha}_ty)+_aD^{\alpha}_t(|_tD^{\alpha}_by|^{p-2}_tD^{\alpha}_b y)\Big)=q(t)|y|^{p-2}y,
 \quad a<t<b,\\[0.3cm]
 y(a)=y(b)=0,
 \end{array}
 \right.
 \end{equation}
 where $_aD^{\alpha}_t$ and $_tD^{\alpha}_b$ denote the left Riemann-Liouville fractional derivative of order $\alpha\in (0,1)$ and  the right Riemann-Liouville fractional derivative of order $\alpha$ respectively. Note that for $\alpha=1,$ problem (\ref{LR-1}) reduces to the $p$-Laplacian problem
\begin{equation}\label{LR-2}
 \left\{\begin{array}{ll}
(|y'|^{p-2}y)'+q(t)|y|^{p-2}y=0,
 \quad a<t<b,~~ p>1,\\
 y(a)=y(b)=0,
 \end{array}
 \right.
 \end{equation}
  The Lyapunov-type inequality for the problem (\ref{LR-1}) is given in  the following theorem.

  \begin{theorem}
  Assume that $0\frac{1}{p}<\alpha<1$ and $q\in L^1(a,b).$ If  problem (\ref{LR-1}) admits a nontrivial solution $y\in AC_a^{\alpha,p}[a,b]\cap AC_b^{\alpha,p}[a,b]\cap C[a,b]$ such that $|y(c)|=\|y\|_{\infty}, ~ c\in (a,b),$ then
  $$\int_a^bq^{+}(s)ds\ge \Bigg(\frac{2(\alpha p-1)}{p-1}\Bigg)^{p-1}\frac{[\Gamma(\alpha)]^p}{\Big((c-a)^{\frac{\alpha p-1}{p-1}}+(b-c)^{\frac{\alpha p-1}{p-1}}\Big)^{p-1}}, $$
  where $q^{+}(t)=\max\{q(t), 0\}$ for $t\in [a,b], $ $AC_a^{\alpha,p}[a,b]=\{y\in L^1(a,b): _aD^{\alpha}_ty\in L^p(a,b)\}$ and $AC_b^{\alpha,p}[a,b]=\{y\in L^1(a,b): _tD^{\alpha}_by\in L^p(a,b)\}.$
  \end{theorem}

\section{Lyapunov inequality for   boundary value problems  with nonsingular Mittag-Leffler
kernel}

In 2017, Abdeljawad \cite{Abd} proved a
Lyapunov-type inequality for the following Riemann-Liouville type fractional boundary value problem
of order $2 < \alpha\le 3$ in terms of Mittag-Leffler kernels:
\begin{equation}\label{p-L-7-x}
 \left\{\begin{array}{ll}
(_a^{ABR}D^{\alpha}y)(t)+q(t)y(t)=0,  \quad a<t<b,~ 2<\alpha\le 3,\\[0.3cm]
\displaystyle y(a)=y(b)=0,
 \end{array}
 \right.
 \end{equation}
 where $_a^{ABR}D^{\alpha}$ denotes the left Riemann-Liouville fractional derivative defined by
 $$(_a^{ABR}D^{\alpha}f)(t)=\frac{B(\alpha)}{1-\alpha}\frac{d}{dt}\int_a^tf(x)E_{\alpha}\Big(-\alpha\frac{(t-x)^{\alpha}}{1-\alpha}\Big)ds,$$
 where $B(\alpha)$ is a normalization function such that $B(0) = B(1) = 1,$ and $E_{\alpha}$ is the generalized Mittag-Leffler function given by $\displaystyle E_{\alpha}(-t^{\alpha})=\sum_{k=0}^{\infty}\frac{(-t)^{\alpha k}}{\Gamma(\alpha k+1)}.$

The integral equation equivalent to  the boundary value problem (\ref{p-L-7-x}) is
 $$y(t)=\int_a^b G(t,s)R(t,y(s))ds,$$
 where
 $$G(t,s)=
\begin{cases}
\displaystyle\frac{(t-a)(b-s)}{b-a}, &a\leq t\leq s\leq b,\\
\displaystyle\frac{(t-a)(b-s)}{b-a}-(t-s), &a\leq s \leq t \leq b,
\end{cases}
$$
and
$$R(t,y(t))=\frac{1-\beta}{B(\beta)}q(t)y(t)+\frac{\beta}{B(\beta)}\big(_aI^{\beta}q(\cdot)y(\cdot)\big)(t), ~~\beta=\alpha-2.$$

The Green's function $G(t, s)$ defined above has the following properties:
\begin{itemize}
\item[(i)] $G(t, s) \ge 0$ for all $a \le t, s \le b;$
\item[(ii)]~ $\max_{t\in [a,b]} G(t, s) = G(s, s)$ for $s \in [a, b];$
\item[(iii)]~ $G(s, s)$ has a unique maximum, given by
$$\max_{s\in [a,b]} G(s, s)=G\Big(\frac{a+b}{2},\frac{a+b}{2}\Big)=\frac{b-a}{4}.$$
\end{itemize}
The Lyapunov inequality for the problem (\ref{p-L-7-x}) is given in the following result.
\begin{theorem}
If the boundary value problem (\ref{p-L-7-x}) has a nontrivial solution, where $q$ is
a real-valued continuous function on $[a, b],$ then
\begin{equation} \label{awa}
\int_a^b\Big[\frac{3-\alpha}{B(\alpha-2)}|q(t)|+\frac{\alpha-2}{B(\alpha-2)}(_aI^{\alpha-2}|q(\cdot)|)(t)\Big]ds>\frac{4}{b-a}.
\end{equation}
\end{theorem}
\begin{remark}
For $\alpha\to 2+,$ notice that  $\displaystyle\frac{3-\alpha}{B(\alpha-2)}|q(t)|+\frac{\alpha-2}{B(\alpha-2)}(_aI^{\alpha-2}|q(\cdot)|)(t) \to |q(t)|$ and hence the inequality (\ref{awa}) reduces to the classical
Lyapunov inequality (\ref{E-5}).
\end{remark}

\section{Lyapunov inequality for discrete fractional  boundary value problems}
Let us begin this section  with the definitions of integral and derivative
of arbitrary order for a function defined on a discrete set. For details, see   \cite{GP}.

The power function is defined by
$$x^{(y)}=\frac{\Gamma(x+1)}{\Gamma(x+1-y)},~ \mbox{for}~ x, x-y\in {\mathbb R}\setminus\{\ldots, -2,-1\}.$$
For $a\in {\mathbb R},$ we define the set ${\mathbb N}_a=\{a, a+1,a+2, \ldots\}.$ Also, we use the notation
$\sigma (s) = s+1$ for the shift operator and $(\Delta f)(t) = f (t +1)- f (t)$ for the forward
difference operator. Notice that $(\Delta^2 f)(t) = (\Delta\Delta f)(t).$

For a function $f : {\mathbb N}_a\to {\mathbb R},$ the discrete fractional sum of order $\alpha\ge 0$ is defined
as
\begin{eqnarray*}
(_a\Delta^0 f)(t)&=&f(t), ~~ t\in {\mathbb N}_a,\\
(_a\Delta^{-\alpha} f)(t)&=&\frac{1}{\Gamma(\alpha)}\sum_{s=a}^{t-\alpha}(t-\sigma(s))^{(\alpha-1)}, f(s), ~t\in {\mathbb N}{a+\alpha},~ \alpha>0.
\end{eqnarray*}

The discrete fractional derivative of order $\alpha\in (1,2]$ is defined by
$$(_a\Delta^{\alpha} f)(t)=(\Delta^2\, _a\Delta^{-(2-\alpha)}f)(t),~~
t\in {\mathbb N}_{a+2-\alpha}.$$

In 2015, Ferreira \cite{r6} studied the following  conjugate boundary value problem
 \begin{equation}\label{D-1}
 \left\{\begin{array}{ll}
\displaystyle  (\Delta^{\alpha}y)(t)+q(t+\alpha-1)y(t+\alpha-1)=0, ~~ t\in [0,b+1]_{\mathbb N_0},\\
\displaystyle y(\alpha-2)=0=y(\alpha+b+1).
 \end{array}
 \right.
 \end{equation}

The function $y$ is a solution of the boundary value problem (\ref{D-1})
if and only if $y$ satisfies the integral equation
\begin{equation}\label{greensol-D}
y(t)=\sum_0^{b+1}G(t,s)q(s+\alpha-1)f(y(s+\alpha-1)),
\end{equation}
where $G(t,s)$ is the Green's function defined by
\begin{equation}\label{green-D}
G(t,s)=\frac{1}{\Gamma(\alpha)}
\begin{cases}
\displaystyle\frac{t^{(\alpha-1)}(\alpha+b-\sigma(s))^{(\alpha-1)}}{(\alpha+b+1)^{(\alpha-1)}}\\
-(t-\sigma(s))^{(\alpha-1)}, &s< t-\alpha+1< b+1,\\[0.3cm]
\displaystyle\frac{t^{(\alpha-1)}(\alpha+b-\sigma(s))^{(\alpha-1)}}{(\alpha+b+1)^{(\alpha-1)}}, &t-\alpha+1\leq s  \leq b+1,
\end{cases}
\end{equation}
and that
$$\max_{s\in [0,b+1]_{\mathbb N_0}}G(s,s)=G\Big(\frac{b}{2}+\alpha-1, \frac{b}{2}\Big),~~\mbox{if}~ b~\mbox{is even},$$
$$\max_{s\in [0,b+1]_{\mathbb N_0}}G(s,s)=G\Big(\frac{b+1}{2}+\alpha-1, \frac{b+1}{2}\Big),~~\mbox{if}~ b~\mbox{is odd},$$
The Lyapunov inequality for the problem (\ref{D-1}) is as follows.
\begin{theorem}\label{t-D-1}
If the discrete fractional boundary value problem (\ref{D-1})
has a nontrivial solution, then
$$\sum_{s=0}^{b+1}|q(s+\alpha-1)|>4\Gamma(\alpha)\frac{\Gamma(b+\alpha+2)\Gamma^2\big(\frac{b}{2}+2\Big)}{(b+2\alpha)(b+2)\Gamma^2\big(\frac{b}{2}+\alpha\Big)\Gamma(b+3)}, ~~\mbox{if}~ b~\mbox{is even};$$
$$\sum_{s=0}^{b+1}|q(s+\alpha-1)|>4\Gamma(\alpha)\frac{\Gamma(b+\alpha+2)\Gamma^2\big(\frac{b+3}{2}\Big)}{\Gamma(b+3)(\Gamma^2\big(\frac{b+1}{2}+\alpha\Big)}, ~~\mbox{if}~ b~\mbox{is odd}.$$
\end{theorem}

 As a simple application, consider the right-focal boundary value problem in Theorem \ref{t-D-1}
with $q = \lambda\in {\mathbb R}.$
Then an eigenvalue of the boundary value problem
 \begin{equation}\label{exam-8}
 \left\{\begin{array}{ll}
(\Delta^{\alpha} y)(t)+\lambda y(t+\alpha-1)=0,  \quad t\in [0, b+1]_{{\mathbb N}_0}, \\[0.2cm]
\displaystyle y(\alpha-2)=0=\Delta y(\alpha+b),
 \end{array}
 \right.
 \end{equation}
must necessarily satisfy the following inequality
$$|\lambda|\ge \frac{1}{\Gamma(\alpha-1)(b+2)^2}.$$

Ferreira \cite{r6} also studied the following right-focal boundary value problem
\begin{equation}\label{D-2}
 \left\{\begin{array}{ll}
\displaystyle  (\Delta^{\alpha}y)(t)+q(t+\alpha-1)y(t+\alpha-1)=0, ~~ t\in [0,b+1]_{\mathbb N_0},\\
\displaystyle y(\alpha-2)=0=\Delta y(\alpha+b).
 \end{array}
 \right.
 \end{equation}
 The function  $y$ is a solution of the boundary value problem (\ref{D-2})
if and only if $y$ satisfies the integral equation
\begin{equation}\label{G-D-2}
y(t)=\int_a^bG(t,s)q(s)y(s)ds,
\end{equation}
where $G(t,s)$ is the Green's function defined by
 {\small \begin{equation}\label{Gr-D-2}
G(t,s)=\frac{1}{\Gamma(\alpha)}
\begin{cases}
\displaystyle\frac{\Gamma(b+3)t^{(\alpha-1)}(\alpha+b-\sigma(s)^{(\alpha-2)}}{\Gamma(\alpha+b+1)}\\
-(t-\sigma(s))^{(\alpha-1)}, &s< t-\alpha+1< b+1,\\[0.3cm]
\displaystyle\frac{\Gamma(b+3)t^{(\alpha-1)}(\alpha+b-\sigma(s)^{(\alpha-2)}}{\Gamma(\alpha+b+1)}, &t-\alpha+1\leq s  \leq b+1,
\end{cases}
\end{equation}}
with
$$\max_{s\in [0,b+1]_{\mathbb N_0}}G(s+\alpha-1,s)=G(b+\alpha, b+1)=\Gamma(\alpha-1)(b+2).$$
The Lyapunov inequality for the problem (\ref{D-2}) is presented as follows.
\begin{theorem}
If the discrete fractional boundary value problem (\ref{D-2})
has a nontrivial solution, then
$$\sum_{s=0}^{b+1}|q(s+\alpha-1)|>\frac{1}{\Gamma(\alpha-1)(b+2)}.$$
\end{theorem}

In 2017, Chidouh and Torres \cite{ChTo-1}, studied the following  conjugate boundary value problem
 \begin{equation}\label{D-1-x}
 \left\{\begin{array}{ll}
\displaystyle  (\Delta^{\alpha}y)(t)+q(t+\alpha-1)f(y(t+\alpha-1))=0, ~~ t\in [0,b+1]_{\mathbb N_0},\\[0.2cm]
\displaystyle y(\alpha-2)=0=y(\alpha+b+1),
 \end{array}
 \right.
 \end{equation}
 where $f\in C({\mathbb R}^+, {\mathbb R}^+)$ is nondecreasing and $q: [\alpha-1, \alpha+b]_{{\mathbb N}_{\alpha-1}}\to {\mathbb R}^+$ is a non-trivial function.

The function $y$ is a solution of the boundary value problem (\ref{D-1-x})
if and only if $y$ satisfies the integral equation
\begin{equation}\label{G-1-x}
y(t)=\sum_0^{b+1}G(t,s)q(s+\alpha-1)f(y(s+\alpha-1)),
\end{equation}
where $G(t,s)$ is the Green's function defined by
 {\small \begin{equation}\label{Gr-1-x}
G(t,s)=\frac{1}{\Gamma(\alpha)}
\begin{cases}
\displaystyle\frac{t^{(\alpha-1)}(\alpha+b-s)^{(\alpha-1)}}{(\alpha+b+1)^{(\alpha-1)}}-(t-s-1)^{(\alpha-1)}, &s< t-\alpha+1< b+1,\\
\displaystyle\frac{t^{(\alpha-1)}(\alpha+b-s)^{(\alpha-1)}}{(\alpha+b+1)^{(\alpha-1)}}, &t-\alpha+1\leq s  \leq b+1.
\end{cases}
\end{equation}}
Moreover, the function $G$ satisfies the following properties:
\begin{itemize}
\item[(i)] $G(t,s)>0$ for all $t\in [\alpha-1, \alpha+b]_{\mathbb N_0}$ and $s\in [1,b+1]_{\mathbb N_1};$
\item[(ii)]~ $\max_{[\alpha-1, \alpha+b]_{\mathbb N_0}}G(t,s)=G(s+\alpha-1,s), ~ s\in [1,b+1]_{\mathbb N_1};$
\item[(iii)]~ $G(s+\alpha-1)$ has a unique maximum given by
$$\max_{ s\in [1,b+1]_{\mathbb N_1}}G(s+\alpha-1)=\begin{cases}
\displaystyle\frac{1}{4} \frac{(b+2\alpha)(b+2)\Gamma^2\Big(\frac{b}{2}+\alpha\Big)\Gamma(b+3)}{\Gamma(\alpha)\Gamma(b+\alpha+2)\Gamma^2\Big(\frac{b}{2}+\alpha\Big)}&if~ b~ \mbox{is even},\\
\displaystyle\frac{1}{\Gamma(\alpha)} \frac{(b+3)\Gamma^2\Big(\frac{b+1}{2}+\alpha\Big)}{\Gamma(b+\alpha+2)\Gamma^2\Big(\frac{b+3}{2}\Big)}&if~ b~ \mbox{is odd}.
\end{cases}
$$
\end{itemize}

The Lyapunov inequality for the problem (\ref{D-1-x}) is expressed as follows.
\begin{theorem}
If the discrete fractional boundary value problem (\ref{D-1-x})
has a nontrivial solution, then
$$\sum_{s=0}^{b+1}|q(s+\alpha-1)|>4\Gamma(\alpha)\frac{\Gamma(b+\alpha+2)\Gamma^2\big(\frac{b}{2}+2\Big)\eta}{(b+2\alpha)(b+2)\Gamma^2\big(\frac{b}{2}+\alpha\Big)\Gamma(b+3)f(\eta)}, ~~\mbox{if}~ b~\mbox{is even},$$
and
$$\sum_{s=0}^{b+1}|q(s+\alpha-1)|>4\Gamma(\alpha)\frac{\Gamma(b+\alpha+2)\Gamma^2\big(\frac{b+3}{2}\Big)\eta}{\Gamma(b+3)(\Gamma^2\big(\frac{b+1}{2}+\alpha\Big)f(\eta)}, ~~\mbox{if}~ b~\mbox{is odd},$$
where
$\eta=\max_{[\alpha-1, \alpha+b]_{\mathbb N_{\alpha-1}}}y(s+\alpha-1).$
\end{theorem}

In 2017, Ghanbari and Gholami \cite{GG} presented the Lyapunov-type inequalities for two special classes of Sturm-Liouville
problems equipped with fractional $\Delta$-difference operators, which are given  in the next two results.
\begin{theorem}
Assume that $p: [a,b]_{{\mathbb N}_0}\to {\mathbb R}^+$  and  $q: [\alpha+a-1,\alpha+b-1]_{{\mathbb N}_{\alpha-1}}\to {\mathbb R}$ are real-valued functions.
If $y$ defined on $[\alpha+a-1,\alpha+b-1]_{{\mathbb N}_{\alpha-1}}$  is a nontrivial
solution to the fractional Sturm-Liouville problem
\begin{equation}\label{exam-19}
\begin{cases}
\Delta^{\alpha}_{b-}(p(t)\Delta^{\alpha}_{a+}y(t))+[q(t+\alpha-1)-\lambda]y(t+\alpha-1)=0, \quad t\in(a,b),\\
y(\alpha+a-1)=0,~ y(\alpha+b)=0,
\end{cases}
\end{equation}
where $\alpha\in (1/2,1)$ and $t=a,a+1,\ldots, b,$ $a,b\in {\mathbb Z},$ $\lambda\in {\mathbb R}$ such that $a\ge 1, b\ge 3,$ then the following Lyapunov type inequality holds:
$$\sum_{s=a}^b\sum_{w=a}^b\Bigg(\frac{|q(w+\alpha-1)-\lambda|}{p(s)}\Bigg)\ge \frac{1}{2}.$$
\end{theorem}
\begin{theorem}
Suppose that   $q: [\alpha+a-1,\alpha+b-1]_{{\mathbb N}_{\alpha-1}}\to {\mathbb R}$ is a real-valued function for $1<\alpha\le 2.$
Assume that  $y$ defined on $[\alpha+a-2,\alpha+b+1]_{{\mathbb N}_{\alpha-2}}$  is a nontrivial
solution to the fractional $\Delta$-difference boundary value  problem:
\begin{equation}\label{exam-19-a}
\begin{cases}
\Delta^{\alpha}_{a+}y(t)+[q(t+\alpha-1)-\lambda]y(t+\alpha-1)=0, \quad t\in(a,b),\\
y(\alpha+a-2)=0,~ y(\alpha+b+1)=0,
\end{cases}
\end{equation}
where   $t=a,a+1,\ldots, b,b+1$ $a,b\in {\mathbb Z},$ $\lambda\in {\mathbb R}$ such that $a\ge 1, b\ge 2,$
then the following Lyapunov-type inequalities hold:
$$\sum_{a}^{b+1}|q(s+\alpha-1)-\lambda|\ge \Gamma(\alpha)\frac{b-a+2}{b-a+2\alpha}\frac{\Gamma(\alpha+b-a+2)\Gamma^2\Big(\frac{b-a}{2}+1\Big)}{\Gamma(b-a+2)\Gamma^2\Big(\frac{b-a}{2}+\alpha\Big)},$$
if $a+b$ is even and
$$\sum_{a}^{b+1}|q(s+\alpha-1)-\lambda|\ge \Gamma(\alpha)\frac{b-a+3}{b-a+2\alpha+1}\frac{\Gamma(\alpha+b-a+2)\Gamma^2\Big(\frac{b-a+1}{2}+1\Big)}{\Gamma(b-a+2)\Gamma^2\Big(\frac{b-a+1}{2}+\alpha\Big)},$$
if $a+b$ is odd.
\end{theorem}

\section{Lyapunov-type inequality for fractional impulsive boundary value problems}
In 2017,   Kayar \cite{Kayar} considered the following impulsive fractional boundary value problem
\begin{equation}\label{exam-14}
 \left\{\begin{array}{ll}
 (^CD^{\alpha} y)(t)+q(t) y(t)=0,  \quad t\ne \tau_i,~~a<t<b,  1<\alpha<2, \\
\displaystyle  \Delta y|_{t=\tau_i}:=y(\tau_i^+)-y(\tau_i^-), ~~ i=1,2,\ldots, p,\\
\displaystyle \Delta y'_{t=\tau_i}=-\frac{\gamma_i}{\beta_i}y(\tau_i^-), ~~ i=1,2,\ldots, p,\\
y(a)=y(b)=0,
 \end{array}
 \right.
 \end{equation}
where $^CD^{\alpha}$ is the Caputo fractional derivative of order $\alpha$ ($1<\alpha\le 2$), $q: PLC[a,b]\to {\mathbb R}$ is a continuous function,
$a=\tau_0<\tau_1<\ldots<\tau_p<\tau_{p+1}=b,$ $PLC[a,b]=\{y: [a,b]\to {\mathbb R}$ is continuous on each interval $(\tau_i, \tau_{i+1}),$ the limits
$y(\tau_i^{\pm})$ exist and $y(\tau_i^-)=y(\tau_i)$ for $i=1,2,\ldots,p\},$ and $PLC^1[a,b]=\{y: [a,b]\to {\mathbb R}, y'\in PLC[a,b]\}.$

$y\in PLC^1([a, b], {\mathbb R})$ is a solution of the boundary value problem (\ref{exam-14})  if and only if y satisfies the following integral equation
$$
y(t)=-\int_a^bG(t,s)q(s)y(s)ds-\sum_{a\le \tau_i<b}H(t,\tau_i)\frac{\gamma_i}{\beta_i}y(\tau_i),
$$
where
 $$G(t,s)=\frac{1}{\Gamma(\alpha)}
\begin{cases}
\displaystyle \frac{a-t}{b-a}(b-s)^{\alpha-1} &a\leq t\leq s\leq b,\\[0.4cm]
\displaystyle \frac{a-t}{b-a}(b-s)^{\alpha-1}-(t-s)^{\alpha-1}, &a\leq s \leq t \leq b,
\end{cases}
$$
and
$$H(t,\tau_i)=\frac{1}{\Gamma(\alpha)}
\begin{cases}
\displaystyle \frac{a-t}{b-a}(b-\tau_i) &a\leq t\leq \tau_i\leq b,\\[0.4cm]
\displaystyle \frac{a-\tau_i}{b-a}(b-t), &a\leq \tau_i \leq t \leq b.
\end{cases}
$$
Furthermore, the functions $G$ and $H$ satisfy the following properties:
\begin{itemize}
\item[(i)] $G(t,s)|\le \frac{1}{\Gamma(\alpha)}\frac{(\alpha-1)^{\alpha-1}}{\alpha^{\alpha}}(b-a)^{\alpha-1}, $ for all $a\le t,s\le b;$
\item[(ii)]~ $H(t,\tau_i)\le 0$ and $|H(t,\tau_i)|\le \frac{b-a}{4}, $ for all $a\le t, \tau_i\le b.$
\end{itemize}
The Lyapunov inequality for the problem (\ref{exam-14}) is the following.
\begin{theorem} (Lyapunov inequality)
If the problem (\ref{exam-14}) has a nontrivial solution $y(t)\ne 0$ on $(a,b),$ then
$$\int_a^b|q(s)|ds+\sum_{a\le \tau_i<b}\Big(\frac{\gamma_i}{\beta_i}\Big)^{+}>\min\Bigg\{\frac{4}{b-a}, \frac{\Gamma(\alpha)\alpha^{\alpha}}{[(\alpha-1)(b-a)]^{\alpha-1}}\Bigg\},$$
where
$\Big(\frac{\gamma_i}{\beta_i}\Big)^{+}=\max\Big\{\frac{\gamma_i}{\beta_i},0\Big\}.$
\end{theorem}

 \section{Lyapunov inequality for  boundary value problems
involving Hilfer fractional  derivative}
A generalization of  both Riemann-Liouville and Caputo derivatives, known as   {\em generalized Riemann-Liouville derivative of order $\alpha\in  (0, 1)$  and  type $\beta\in  [0, 1]$},   was proposed  by R. Hilfer  in \cite{Hil}.   Such a derivative interpolates between the Riemann-Liouville and Caputo derivative in some sense. For properties and applications of the Hilfer derivative, see \cite{Hil-1, HKT} and references cited therein.

\begin{definition}\label{gus3}
	The  generalized Riemann-Liouville fractional derivative or Hilfer fractional derivative of order $\alpha$ and parameter $\beta$ for  a function $u$ is defined as
	\begin{eqnarray*}
		^H D^{\alpha,\beta} u(t)=I^{\beta (n-\alpha)}D^nI^{(1-\beta)(n-\alpha)}u(t),
	\end{eqnarray*}
	where $n-1<\alpha < n ,\hspace{0.2cm} 0\leq \beta \leq 1, \hspace{0.2cm }t>a ,\hspace{0.2cm} D=\displaystyle\frac{d}{dt}$.
\end{definition}

\begin{remark}
	The Hilfer fractional derivative corresponds to the Riemann-Liouville fractional derivative for $\beta =0,$ that is,
	$	^H D^{\alpha,0}u(t)=D^n I^{n-\alpha}u(t),$
	while  it corresponds to the Caputo fractional derivative for $\beta=1$ given by
	$
	^H D^{\alpha,1}u(t)=I^{n-\alpha} D^n u(t).
	$
\end{remark}

In 2016, Pathak \cite{Path} studied Lyapunov-type inequalities for fractional boundary value problems involving Hilfer fractional derivative and Dirichlet,  mixed Dirichlet and Neumann boundary conditions.

 Let us first consider the Dirichlet boundary value problem given by
\begin{equation}\label{exam-17}
\begin{cases}
\displaystyle (D^{\alpha,\beta}y)(t)+g(t)y(t)=0, \quad t\in(a,b),~ 1<\alpha\le 2, ~ 0\le \beta\le 1,\\[0.3cm]
y(a)=y(b)=0,
\end{cases}
\end{equation}
which is equivalent to the integral equation
$$y(t)=\int_a^bG(t,s)q(s)y(s)ds,$$
where
$$G(t,s)=\frac{1}{\Gamma(\alpha)}
\begin{cases}
\displaystyle\Bigg(\frac{t-a}{b-a}\Bigg)^{1-(2-\alpha)(1-\beta)}(b-s)^{\alpha-1}, &a\leq t\leq s\leq b,\\[0.4cm]
\displaystyle \Bigg(\frac{t-a}{b-a}\Bigg)^{1-(2-\alpha)(1-\beta)}(b-s)^{\alpha-1}-(t-s)^{\alpha-1}, &a\leq s \leq t \leq b,
\end{cases}
$$ is the Green's function satisfying the property:
$$|G(t,s)|\le \frac{(b-a)^{\alpha-1}[\alpha-1+\beta(2-\alpha)]^{\alpha-1+\beta(2-\alpha)}[\alpha-1]^{\alpha-1}}{\Gamma(\alpha)[\alpha-(2-\alpha)(1-\beta)]^{\alpha-(2-\alpha)(1-\beta)}}, ~ (t,s)\in [a,b]\times [a,b].$$

\begin{theorem} (Lyapunov-type inequality)
If a nontrivial continuous solution of the problem (\ref{exam-17}) exists, then
$$\int_a^b|q(s)|ds\ge\frac{\Gamma(\alpha)[\alpha-(2-\alpha)(1-\beta)]^{\alpha-(2-\alpha)(1-\beta)}}{(b-a)^{\alpha-1}[\alpha-1+\beta(2-\alpha)]^{\alpha-1+\beta(2-\alpha)}(\alpha-1)^{\alpha-1}}.$$
\end{theorem}

Next we consider a fractional boundary value problems involving Hilfer fractional derivative and   mixed Dirichlet and Neumann boundary conditions:
\begin{equation}\label{exam-17a}
\begin{cases}
\displaystyle (D^{\alpha,\beta}y)(t)+q(t)y(t)=0, \quad t\in(a,b),~ 1<\alpha\le 2, ~ 0\le \beta\le 1,\\[0.2cm]
y(a)=y'(b)=0,
\end{cases}
\end{equation}
which is equivalent to the integral equation
$$y(t)=\int_a^bG(t,s)q(s)y(s)ds,$$
where $\displaystyle G(t,s)=\frac{H(t,s)}{\Gamma(\alpha)(b-s)^{2-\alpha}}$
and
$$H(t,s)=
\begin{cases}
\displaystyle\frac{(\alpha-1)(t-a)^{1-(2-\alpha)(1-\beta)}(b-a)^{(2-\alpha)(1-\beta)}}{1-(2-\alpha)(1-\beta)}, &a\leq t\leq s\leq b,\\[0.4cm]
\displaystyle\frac{(\alpha-1)(t-a)^{1-(2-\alpha)(1-\beta)}(b-a)^{(2-\alpha)(1-\beta)}}{1-(2-\alpha)(1-\beta)}\\
-(t-s)^{\alpha-1}(b-s)^{2-\alpha}, &a\leq s \leq t \leq b.
\end{cases}
$$
The function $H$   satisfies the following property:
$$|H(t,s)|\le\frac{b-a}{\alpha-1+\beta(2-\alpha)}\max\{\alpha-1,\beta(2-\alpha)\}, ~ (t,s)\in [a,b]\times [a,b].$$

\begin{theorem} (Lyapunov-type inequality)
If a nontrivial continuous solution of the problem (\ref{exam-17a}) exists, then
$$\int_a^b(b-s)^{\alpha-2}|q(s)|ds\ge\frac{\Gamma(\alpha)(\alpha-1+\beta(2-\alpha))}{(b-a)\max\{\alpha-1,\beta(2-\alpha)\}}.$$
\end{theorem}

Now we establish a Lyapunov-type inequality for another  fractional boundary value problems with Hilfer fractional derivative and   mixed Dirichlet and Neumann boundary conditions:
\begin{equation}\label{exam-17b}
\begin{cases}
\displaystyle (D^{\alpha,\beta}y)(t)+q(t)y(t)=0, \quad t\in(a,b),~ 1<\alpha\le 2, ~ 0\le \beta\le 1,\\[0.2cm]
y(a)=y'(a)=y'(b)=0,
\end{cases}
\end{equation}
which can be transformed to the integral equation: $$y(t)=\int_a^bG(t,s)q(s)y(s)ds,$$
where
$$G(t,s)=\frac{1}{\Gamma(\alpha)}
\begin{cases}
\displaystyle\frac{(\alpha-1)(t-a)^{2-(3-\alpha)(1-\beta)}(b-s)^{\alpha-2}}{(b-a)^{1-(3-\alpha)(1-\beta)}[2-(3-\alpha)(1-\beta)}, &a\leq t\leq s\leq b,\\[0.4cm]
\displaystyle\frac{(\alpha-1)(t-a)^{2-(3-\alpha)(1-\beta)}(b-s)^{\alpha-2}}{(b-a)^{1-(3-\alpha)(1-\beta)}[2-(3-\alpha)(1-\beta)}\\
-(t-s)^{\alpha-1}, &a\leq s \leq t \leq b,
\end{cases}
$$ is the Green's function  satisfying the property:
$$|G(t,s)|\le\frac{2(b-a)^{\alpha-1}(\alpha-2)^{\alpha-2}}{\Gamma(\alpha)[2-(3-\alpha)(1-\beta)]^{\alpha-1}}, ~ (t,s)\in [a,b]\times [a,b].$$

\begin{theorem}(Lyapunov-type inequality)
If a nontrivial continuous solution of the problem (\ref{exam-17b}) exists, then
$$\int_a^b |q(s)|ds\ge\frac{\Gamma(\alpha)[2-(3-\alpha)(1-\beta)]^{\alpha-1}}{(b-a)^{\alpha-1}(\alpha-2)^{\alpha-2}}.$$
\end{theorem}
\vskip 0.3cm
Finally we consider the following fractional boundary value problem with Hilfer fractional derivative and a mixed set of fractional Dirichlet, Neumann and fractional Neumann boundary conditions
\begin{equation}\label{HIL-4}
 \left\{\begin{array}{ll}
(D^{\alpha,\beta}y)(t)+q(t)y(t)=0,
 \quad a<t<b,~~ 2<\alpha\le 3, ~ 0\le \beta\le 1,\\[0.3cm]
 \displaystyle \Big(I^{(3-\alpha)(1-\beta)}y\Big)(a)=0, ~ y'(b)=0,~ \frac{d^2}{dt^2}\Big(I^{(3-\alpha)(1-\beta)}y\Big)(a)=0.
 \end{array}
 \right.
 \end{equation}
The integral equation equivalent to the Problem (\ref{HIL-4}) is
$$y(t)=\int_a^bG(t,s)q(s)y(s)ds,$$
where
$$G(t,s)=\frac{1}{\Gamma(\alpha)}
\begin{cases}
\displaystyle \frac{(\alpha-1)(t-a)^{1-(3-\alpha)(1-\beta)}(b-s)^{\alpha-2}}{(b-a)^{-(3-\alpha)(1-\beta)}[1-(3-\alpha)(1-\beta)]}, &a\leq t\leq s\leq b,\\[0.3cm]
\displaystyle \frac{(\alpha-1)(t-a)^{1-(3-\alpha)(1-\beta)}(b-s)^{\alpha-2}}{(b-a)^{-(3-\alpha)(1-\beta)}[1-(3-\alpha)(1-\beta)]}\\
-(t-s)^{\alpha-1}, &a\leq s \leq t \leq b,
\end{cases}
$$ such that
$$|G(t,s)|\le \frac{(\alpha-2)^{\alpha-2}(b-a)^{\alpha-1}}{\Gamma(\alpha)[1-(3-\alpha)(1-\beta)]^{\alpha-1}},~~ (t,s)\in [a,b]\times [a,b].$$

\begin{theorem} (Lyapunov inequality)
If a nontrivial continuous  solution of the problem (\ref{HIL-4}) exists, then
$$\int_a^b|q(s)|ds\ge \frac{\Gamma(\alpha)[1-(3-\alpha)(1-\beta)]^{\alpha-1}}{(b-a)^{\alpha-1}(\alpha-2)^{\alpha-2}}.$$
\end{theorem}

In 2017,  Kirane and Torebek   \cite{KiTo}  obtained   Lyapunov-type inequalities for the following {\em  fractional
 boundary value problem}
  \begin{equation}\label{Jan-1}
 \left\{\begin{array}{ll}
D^{\alpha,\gamma}_{a} y(t)+q(t)f(y(t))=0, \quad a<t<b, ~ 1<\alpha\le \gamma<2,\\[0.2cm]
\displaystyle y(a)=y(b)=0,
 \end{array}
 \right.
 \end{equation}
 where   $D^{\alpha,\gamma}_a$ is  a generalized Hilfer  fractional derivative of  order $\alpha\in {\mathbb R}~(m-1<\alpha<m, m\in {\mathbb N})$ and type $\gamma$,  defined as
 $$D^{\alpha,\gamma}_af(t)=I^{\gamma-\alpha}_a\frac{d^m}{dt^m}I^{m-\gamma}_af(t),$$
  and  $q: [a,b]\to {\mathbb R}$ is a nontrivial Lebesgue  integrable function.

 The integral representation for the solution of the boundary value problem (\ref{Jan-1}) is
\begin{equation*}
y(t)=\int_a^bG(t,s)q(s)f(y(s))ds,
\end{equation*}
where $G(t,s)$ is the Green's function given by
$$G(t,s)=
\begin{cases}
\displaystyle \Big(\frac{t-a}{b-a}\Big)^{\gamma-1}\frac{(b-s)^{\alpha-1}}{\Gamma(\alpha)}-\frac{(t-s)^{\alpha-1}}{\Gamma(\alpha)}, &a\leq s\leq t\leq b,\\[0.4cm]
\displaystyle  \Big(\frac{t-a}{b-a}\Big)^{\gamma-1}\frac{(b-s)^{\alpha-1}}{\Gamma(\alpha)}, &a\leq t \leq s \leq b.
\end{cases}
$$

Further, the above Green's function $G(t,s)$  satisfies the following properties:
\begin{itemize}
\item[(i)]~ $G(t,s)\ge 0$ for $a\le t,s\le b;$
\item[(ii)]~ $\max_{a\le t\le b}G(t,s)=G(s,s), ~ s\in [a,b];$
\item[(iii)]~$G(s,s)$ has a unique maximum, given gy
$$\max_{a\le s\le b}G(s,s)=\frac{(\alpha-1)^{\alpha-1}}{(\gamma+\alpha-2)^{\gamma+\alpha-2}}\frac{((\gamma-1)b-(\alpha-1)a)^{\gamma-1}}{\Gamma(\alpha)(b-a)^{\gamma-\alpha}}.$$
\end{itemize}

They obtained the following Lyapunov-type inequalities.
\begin{theorem}\label{t-Jan-1}
If the fractional boundary value problem \eqref{Jan-1}
 has a nontrivial solution for a real valued continuous function $q$, then
 $$\int_a^b|q(s)|ds>\frac{(\gamma+\alpha-2)^{\gamma+\alpha-2}}{(\alpha-1)^{\alpha-1}}\frac{\Gamma(\alpha)(b-a)^{\gamma-\alpha}}{((\gamma-1)b-(\alpha-1)a)^{\gamma-1}}.$$
\end{theorem}

\begin{theorem}\label{th-2-Jan-1}
Let $q:[a,b]\to {\mathbb R}$ be a real notrivial Lebesgue integrable function and $f\in C({\mathbb R}_{+},{\mathbb R}_{+})$ be a concave and nondecreasing function. If there exists a nontrivial solution $y$ for  the problem (\ref{Jan-1}),  then
$$\int_a^b|q(s)|ds>\frac{(\gamma+\alpha-2)^{\gamma+\alpha-2}}{(\alpha-1)^{\alpha-1}}\frac{\Gamma(\alpha)(b-a)^{\gamma-\alpha}}{((\gamma-1)b-(\alpha-1)a)^{\gamma-1}}\frac{\omega}{f(\omega)},$$
where $\omega=\max_{t\in [a,b]}y(t).$
\end{theorem}

\begin{theorem} {\rm (Hartman-Wintner type inequality)}
Let the functions $q$ and $f$ satisfy the conditions of Theorem \ref{th-2-Jan-1}. Suppose that the fractional boundary value problem  (\ref{Jan-1}) has a nontrivial solution. Then
$$\int_a^b(s-a)^{\gamma-1}(b-s)^{\alpha-1}q^+(s)ds>\frac{\|y\|}{f(\|y\|)}\Gamma(\alpha)(b-a)^{\gamma-1}.$$
\end{theorem}

\begin{corollary}
If $f(y)=y$ (linear case) and $q\in L^1([a,b], {\mathbb R}_{+}), $ then
$$\int_a^b(s-a)^{\gamma-1}(b-s)^{\alpha-1}q^+(s)ds>\Gamma(\alpha)(b-a)^{\gamma-1}.$$
\end{corollary}

\section{Lyapunov-type inequality with the Katugampola fractional derivative}

In 2018,   Lupinska and  Odzijewicz   \cite{LuOd}  obtained   a Lyapunov-type inequality for the following {\em  fractional
 boundary value problem}
  \begin{equation}\label{Feb-1}
 \left\{\begin{array}{ll}
D^{\alpha,\rho}_{a+} y(t)+q(t)y(t)=0, \quad a<t<b, ~ \alpha>0, \rho>0,\\[0.2cm]
\displaystyle y(a)=y(b)=0,
 \end{array}
 \right.
 \end{equation}
 where   $D^{\alpha,\rho}_{a+}$ is  the Katugampola  fractional derivative of  order $\alpha$,  defined as
 $$D^{\alpha,\gamma}_{a+}f(t)=\Big(t^{1-\alpha}\frac{d}{dt}\Big)^nI^{n-\alpha}_{a+}f(t),$$
 for $t\in (a,b),$ $n=[\alpha]+1,$  $0<a<t<b\le \infty$  and  $q: [a,b]\to {\mathbb R}$ is a continuous  function. Here $I^{\alpha,\rho}_{a+}$ is the Katugampola  fractional integral defined by
 $$I^{\alpha,\rho}_{a+}f(t)=\frac{\rho^{1-\alpha}}{\Gamma(\alpha)}\int_a^t\frac{s^{\rho-1}}{(t^{\rho}-s^{\rho})^{1-\alpha}}f(s)ds.$$

 The integral representation for the solution of the boundary value problem (\ref{Feb-1}) is
\begin{equation*}
y(t)=\int_a^bG(t,s)q(s)y(s)ds,
\end{equation*}
where $G(t,s)$ is the Green's function given by
$$G(t,s)=\frac{\rho^{1-\alpha}}{\Gamma(\alpha)}
\begin{cases}
\displaystyle \frac{s^{\rho-1}}{(b^{\rho}-s^{\rho})^{1-\alpha}}\Big(\frac{t^{\rho}-a^{\rho}}{b^{\rho}-a^{\rho}}\Big)^{\alpha-1}, &a\leq t\leq s\leq b,\\[0.4cm]
\displaystyle  \frac{s^{\rho-1}}{(b^{\rho}-s^{\rho})^{1-\alpha}}\Big(\frac{t^{\rho}-a^{\rho}}{b^{\rho}-a^{\rho}}\Big)^{\alpha-1}-\frac{s^{\rho-1}}{(t^{\rho}-s^{\rho})^{1-\alpha}}, &a\leq s \leq t \leq b,
\end{cases}
$$
which satisfies the following properties:
\begin{itemize}
\item[(i)]~ $G(t,s)\ge 0$ for $a\le t,s\le b;$
\item[(ii)]~ $\displaystyle \max_{a\le t\le b}G(t,s)=G(s,s)\le  \frac{\max\{a^{\rho-1},b^{\rho-1}\}}{\Gamma(\alpha)}\Big(\frac{b^{\rho}-a^{\rho}}{4\rho}\Big)^{\alpha-1}, ~ s\in [a,b].$
\end{itemize}

They obtained the following Lyapunov-type inequality.
\begin{theorem}\label{t-Feb-1}
If the fractional boundary value problem \eqref{Feb-1}
 has a nontrivial solution for a real valued continuous function $q$, then
 $$\int_a^b|q(s)|ds>\frac{\Gamma(\alpha)}{\max\{a^{\rho-1},b^{\rho-1}\}}\Big(\frac{4\rho}{b^{\rho}-a^{\rho}}\Big)^{\alpha-1}.$$
\end{theorem}

\begin{remark}
In the special case when $\rho=1$ in Theorem \ref{t-Feb-1}, we get the following result
$$\int_a^b|q(s)|ds\ge \Gamma(\alpha)\Big(\frac{4}{b-a}\Big)^{\alpha-1},$$
which is Theorem \ref{L=0}, while taking $\rho\to 0^+$ in Theorem \ref{t-Feb-1},  we have the Lyapunov's  type inequality for the Hadamard fractional derivative:
$$\int_a^b|q(s)|ds\ge \alpha\Gamma(\alpha)\Big(\frac{\log(b/a)}{4}\Big)^{1-\alpha}.$$
\end{remark}

\section{Lyapunov inequality for a boundary value problem
involving the conformable derivative}

Recently, Khalil et al.  \cite{KHYS} introduced a new derivative, which appears in the form of a limit  like the classical derivative and is known as the  conformable derivative. Later, this new local derivative was improved by
Abdeljawad \cite{Abdel}. The importance of the conformable derivative is that it has  properties similar to  the ones of the classical derivative. However, the conformable derivative does not satisfy the index law \cite{Kat, Ort} and the zero order derivative property, that is, the zero order derivative of a differentiable function does not return to the function itself.

In 2017, Khaldi {\em et al.} \cite{KhLa} obtained a Lyapunov-type inequality for the following boundary value problem involving the conformable derivative of order  $1<\alpha<2$ and Dirichlet boundary conditions:
\begin{equation}\label{exam-18}
\begin{cases}
\displaystyle T_{\alpha}y(t)+q(t)y(t)=0, \quad t\in(a,b),\\
y(a)=y'(b)=0,
\end{cases}
\end{equation}
where $T_{\alpha}$ denotes the conformable derivative of order $\alpha$ and $q: [a,b]\to {\mathbb R}$ is  a real continuous function.

The conformable derivative of order $0 <\alpha < 1$ for a function $g : [a,\infty)\to {\mathbb R}$ is defined by
$$T_{\alpha}g(t)=\lim_{\varepsilon\to 0}\frac{g\Big(t+\varepsilon(t-a)^{1-\alpha}\Big)-g(t)}{\varepsilon},~~ t>a.$$
If $T_{\alpha}g(t)$ exists on $(a,b), b>a$ and $\lim_{t\to a+}T_{\alpha}g(t)$ exists, then we define $T_{\alpha}g(a)=\lim_{t\to a+}T_{\alpha}g(t).$

The conformable derivative of order $n <\alpha <n+ 1$ of a function $g : [a,\infty)\to {\mathbb R},$ when $g^{(n)}$ exists,  is defined as
$$T_{\alpha}g(t)=T_{\beta}g^{(n)}(t),$$
where $\beta=\alpha-n\in (0,1).$

The solution $y$ of the problem (\ref{exam-18}) can be written as
$$y(t)=\int_a^bG(t,s)q(s)y(s)ds,$$
where
$$G(t,s)=\frac{1}{b-a}
\begin{cases}
(b-s)(t-a), &a\leq t\leq s\leq b,\\
-(b-a)(t-s)+(b-s)(t-a), &a\leq s \leq t \leq b,
\end{cases}
$$ is
the Green's function, which is nonnegative, continuous and satisfies the property:
$$0\le G(t,s)\le b-a, ~ \mbox{for all}~ t,s\in [a,b].$$
\begin{theorem} (Lyapunov inequality)
Let $q\in C([a,b], {\mathbb R}).$ If the boundary value problem (\ref{exam-18}) has a solution $y\in  AC^2 ([a,b], {\mathbb R})$ such that $y(t)\ne 0$ a.e.
on $(a,b),$ then
$$\int_a^b|q(s)|(s-a)^{\alpha-2}ds\ge \frac{4}{b-a}.$$
\end{theorem}

In 2017,  Abdeljawad {\em et al.} \cite{AAJ}  obtained Lyapunov-type inequality for a Dirichlet boundary value problem involving conformable derivative of order  $1<\alpha<2$:
\begin{equation}\label{Oct-4}
\begin{cases}
\displaystyle T_{\alpha}y(t)+q(t)y(t)=0, \quad t\in(a,b),\\
y(a)=y(b)=0,
\end{cases}
\end{equation}
where $T_{\alpha}$ denotes the conformable derivative of order $\alpha$ and $q: [a,b]\to {\mathbb R}$ is  a real continuous function.

 The solution for the boundary value problem (\ref{Oct-4}) is
\begin{equation*}
y(t)=\int_a^bG(t,s)q(s)y(s)ds,
\end{equation*}
where $G(t,s)$ is the Green's function given by
$$G(t,s)=
\begin{cases}
\displaystyle \frac{(t-a)(b-s)}{b-a}\cdot (s-a)^{\alpha-2}, &a\leq t\leq s\leq b,\\[0.3cm]
\displaystyle  \Big(\frac{(t-a)(b-s)}{b-a}-(t-s)\Big)\cdot (s-a)^{\alpha-2}, &a\leq s \leq t \leq b,
\end{cases}
$$
which satisfies the  properties:
\begin{itemize}
\item[(i)] ~$G(t,s)\ge 0$ for all $a\le t,s\le b;$
\item[(ii)] ~$\max_{t\in [a,b]}G(t,s)=G(s,s)$ for $s\in [a,b];$
\item[(iii)]~ $G(t,s)$  has a unique maximum, given by
$$\max_{s\in [a,b]}G(s,s)=G\Big(\frac{a+(\alpha-1)b}{\alpha},\frac{a+(\alpha-1)b}{\alpha}\Big)=\frac{(b-a)^{\alpha-1}(\alpha-1)^{\alpha-1}}{\alpha^{\alpha}}.$$
\end{itemize}

The Lyapunov inequality for the problem (\ref{Oct-4}) is given in the following result.
\begin{theorem}\label{t-Oct-4}
If the problem (\ref{Oct-4})  has a nontrivial   solution, where $q$  is a real valued continuous function on $[a,b],$ then
$$\int_0^1|q(s)|ds>\frac{\alpha^{\alpha}}{(b-a)^{\alpha-1}(\alpha-1)^{\alpha-1}}.$$
\end{theorem}

\end{document}